\documentclass[AMA,STIX1COL]{WileyNJD-v2}

\usepackage{lipsum}
\usepackage{amsfonts}
\usepackage{graphicx}
\usepackage{epstopdf}

\usepackage[algo2e,ruled,linesnumbered,algosection]{algorithm2e} % User defined

\usepackage{bm}
\usepackage{subfigure}
\usepackage{xcolor}

\hypersetup{
    colorlinks=false,
    pdfborder={0 0 0},
}

\ifpdf
  \DeclareGraphicsExtensions{.eps,.pdf,.png,.jpg}
\else
  \DeclareGraphicsExtensions{.eps}
\fi

\usepackage{enumitem}
\setlist[enumerate]{leftmargin=.5in}
\setlist[itemize]{leftmargin=.5in}

\ifpdf
\hypersetup{
  pdftitle={Deep convolutional recurrent autoencoders for learning low-dimensional feature dynamics of fluid systems},
  pdfauthor={F. J. Gonzalez and M. Balajewicz}
}
\fi

\articletype{Research Article}

\received{}
\revised{}
\accepted{}

\begin{document}

\title{Deep convolutional recurrent autoencoders for learning low-dimensional feature dynamics of fluid systems}

\author[1]{Francisco J. Gonzalez*}
\author[2]{Maciej Balajewicz}

\authormark{Gonzalez and Balajewicz}

\address{\orgdiv{Department of Aerospace Engineering}, \orgname{University of Illinois at Urbana-Champaign}, \orgaddress{\city{Urbana}, \state{Illinois}, \country{USA}}}

\corres{*Francisco J. Gonzalez, \orgdiv{Department of Aerospace Engineering}, \orgname{University of Illinois at Urbana-Champaign}, \orgaddress{\city{Urbana}, \state{Illinois}, \country{USA}}  \email{fjgonza2@illinois.edu}}

\abstract[Summary]{
  Model reduction of high-dimensional dynamical systems alleviates computational
  burdens faced in various tasks from design optimization to model predictive
  control. One popular model reduction approach is based on projecting the
  governing equations onto a subspace spanned by basis functions obtained from
  the compression of a dataset of solution snapshots. However, this method is
  intrusive since the projection requires access to the system operators.
  Further, some systems may require special treatment of nonlinearities to
  ensure computational efficiency or additional modeling to preserve stability.
  In this work we propose a deep learning-based strategy for nonlinear model
  reduction that is inspired by projection-based model reduction where the idea
  is to identify some optimal low-dimensional representation and evolve it in
  time. Our approach constructs a modular model consisting of a deep
  convolutional autoencoder and a modified LSTM network. The deep convolutional
  autoencoder returns a low-dimensional representation in terms of coordinates
  on some expressive nonlinear data-supporting manifold. The dynamics on this
  manifold are then modeled by the modified LSTM network in a computationally
  efficient manner. An offline unsupervised training strategy that exploits the
  model modularity is also developed. We demonstrate our model on three
  illustrative examples each highlighting the model's performance in prediction
  tasks for fluid systems with large parameter-variations and its stability in
  long-term prediction.
}


\maketitle

\section{Introduction}
\label{sec:introduction}
Dynamical systems are used to describe the rich and complex evolution of many
real-world processes. Modeling the dynamics of physical, engineering, and
biological systems is thus of great importance in their analysis, design, and
control. Many fields, such as the physical sciences, are in the fortunate
position of having a first-principles models that describes the evolution of
certain systems with near-perfect accuracy (e.g., the Navier-Stokes equations in
fluid mechanics, or Schr\"odingers equations in quantum mechanics). Although, in
principle it is possible to numerically solve these equations through direct
numerical simulations (DNS), this often yields systems of equations with
millions or billions of degrees of freedom. Even with recent advances in
computational power and memory capacity, solving these high-fidelity models
(HFMs) is still computationally intractable for \textit{multi-query} and
\textit{time-critical} applications such as design optimization, uncertainty
quantification, and model predictive control. Model reduction aims to alleviate
this burden by constructing reduced order models (ROMs) that capture the
large-scale system behavior while retaining physical fidelity.

Some fields, however, such as finance and neuroscience, lack governing laws
thereby restricting the applicability of principled strategies for constructing
low-order models. In recent years, the rise in machine learning and big data
have driven a shift in the way complex spatiotemporal systems are modeled
~\cite{Raissi2017,Brunton2017,Bongard2007,Schaeffer2017,Tran2017}. The abundance
of data have facilitated the construction of so called data-driven models of
systems lacking high-fidelity governing laws. In areas where HFMs do exist,
data-driven methods have become an increasingly popular approach to tackle
previously challenging problems wherein solutions are \textit{learned} from
physical or numerical data ~\cite{Raissi2016a,Brunton2016,San2017}.

In model reduction, machine learning strategies have recently been applied to
many remaining challenges, including learning stabilizing closure terms in
unstable POD-Galerkin models ~\cite{San2017, Benosman2017}, and data-driven model
identification for truncated generalized POD coordinates ~\cite{Wang2017a,
Wang2018, Kani2017a}. A more recent approach involved learning a set of
observable functions spanning a Koopman invariant subspace from which low-order
linear dynamics of nonlinear systems are modeled ~\cite{Otto2017a}. These
approaches constitute just a small portion of the outstanding challenges in
which machine learning can aid in modeling low-dimensional dynamics of complex
systems.

In this work we make progress to this end by proposing a method that uses a
completely data-driven approach to identify and evolve a low-dimensional
representation of a spatiotemporal system. In particular, we employ a deep
convolutional autoencoder to learn an optimal low-dimensional representation of
the full state of the system in the form of a feature vector, or coordinates of
some low-dimensional nonlinear manifold. The dynamics on this manifold are then
learned using a recurrent neural network trained jointly with the autoencoder in
an end-to-end fashion using a set of finite-time trajectories of the system.

\subsection{Reduced order and surrogate modeling}
Model order reduction is part of a broader family of surrogate modeling
strategies that attempt to reduce the computational burden of solving HFMs by
instead solving approximate, low-complexity models. Surrogate models can be
broadly classified into three groups: 1) data-fit models, 2) hierarchical
models, and 3) projection-based model reduction ~\cite{Benner2015, Carlberg2011}.
Data-fit models use simulation or experimental data to fit an input-output map
as a function of system parameters. Some examples include models based on
Gaussian processes ~\cite{Raissi2017, Parish2016}, and feed-forward neural
networks ~\cite{San2017}. Hierarchical or low-fidelity models substitute the HFM
with a lower-fidelity physics-based model that makes simplifying physics
assumptions (e.g., ignoring viscous effects in a fluid flow), use coarser
computational grids, or relaxes solver-tolerances.

In contrast to the first two surrogate modeling approaches, projection-based
model reduction works by directly exploiting the low-dimensional behavior
inherent in many high-dimensional dynamical systems. These methods approximate
the state of the system by an affine trial subspace, and project the HFM onto a
test subspace resulting in a square system of dimension much smaller than the
original high dimension. Over the years, a large variety of empirically-based
approaches for generating the trial and test subspaces have been developed,
including proper orthogonal decomposition (POD) ~\cite{Lumley1970,Holmes1996},
Krylov subspace methods ~\cite{Bai2002}, and dynamic mode decomposition
~\cite{Schmid2010}. Despite the successes of projection-based model reduction,
there exist a number of issues limiting the applicability of these methods.

One issue is that although the projection step effectively constrains the HFM to
a lower dimensional subspace, this does not necessarily provide computational
efficiency for general nonlinear models. Systems with generic, nonpolynomial
nonlinearities or time-varying parameters require an additional layer of
approximation, or \textit{hyper-reduction}, to gain a computaional speed
up.\footnotemark \footnotetext{For linear time-invariant systems, or systems
with polynomial nonlinearities all projection coefficient can be precomputed
offline.} Some approaches for the treatment of nonlinearities include ROMs based
on discrete empirical interpolation (DEIM) ~\cite{Chaturantabut2010}, or
Gauss-Newton with approximated tensors (GNAT) ~\cite{Carlberg2013}. Other methods
employ patchwork of local state space approximations at multiple locations
including piecewise trajectory linearization (TPWL) ~\cite{Rewienski2006} and
ROMs based on trajectory piecewise quadratic (TPWQ) approximations
~\cite{Trehan2016}. A second well known issue, particularly when dealing with
high-Reynolds number fluid flows, is that of stability. POD-based ROMs are
biased towards large energy-producing scales and are not endowed with the small
energy-dissipating scales that maybe dynamically significant
~\cite{Balajewicz2013}. Moreover, projection-based model reduction has the major
disadvantage of being intrusive, requiring  access to the system operators
during the projection step. Thus, while optimal, say POD-based, approximations
can be made of any dataset, the projection step is still limited to systems with
existing governing laws.

\subsection{Contributions and outline}
In this work we develop a deep learning-based nonlinear model reduction strategy
which is completely data-driven. This method employs a deep convolutional
autoencoder to learn an optimal low-dimensional representation of each solution
snapshot and later evolves this representation in time using a type of recurrent
neural network (RNN) called a long short-term memory (LSTM) network. This work
has important similarities to previous work using RNNs to evolve reduced order
models ~\cite{Wang2017a, Kani2017a} and work that employs autoencoders for
dimensionality reduction ~\cite{Hinton2006, Wang2016b, Hartman2017a, Otto2017a}.

Although previous work regarding neural-network based reduced order models
has shown great promise, a number of significant issues remain. Notably, while
deep fully-connected autoencoders, such as the ones employed in ~\cite{Otto2017a,
Wang2016b, Hartman2017a} work well for small systems with a few thousand degrees
of freedom, this approach alone is not scalable as input data increases to
DNS-level sizes (e.g., $10^6-10^9$ degrees of freedom (dof)). For example, an
autoencoder just a single layer reducing input data from $10^6$ dof to $100$
will require training well over $10^8$ parameters, a feat that quickly becomes
computationally intractable as the autoencoder increases in depth.

To avoid this \textit{curse of dimensionality}, we instead propose a
convolutional recurrent autoencoding model that differs significantly from
existing autoencoder-based model reduction approaches in two main ways:
\begin{itemize}
  \item[(i)]{We propose an autoencoding method that exploits local,
  location-invariant correlations present in physical data through the use of
  convolutional neural networks. That is, rather of applying a fully-connected
  autoencoder to the high-dimensional input data we instead apply it to a
  vectorized feature map produced by a convolutional encoder, and similarly the
  reverse is done for reconstruction. The result is the identification of an
  expressive low-dimensional manifold obtained at a much lower cost while
  offering specific advantages over both traditional POD-based ROMs and
  fully-connected autoencoders.}
  \item[(ii)]{We propose a modified LSTM network to model the evolution of
  low-dimensional data representations on this manifold that avoids costly state
  reconstructions at every step. In doing this, we ensure that the evaluation of
  new steps scales only with the size of the low-dimensional representation and
  not with the size of the full dimensional data, which may be large for some
  problems.}
\end{itemize}
Taken together this end-to-end approach both identifies an optimal
low-dimensional representation of a high-dimensional spatiotemporal dataset and
models its dynamics on the underlying data-supporting manifold. Additionally a
two-step unsupervised training strategy is developed that exploits the modularity of the
convolution recurrent autoencoder model.

The paper is organized as follows. ~\autoref{sec:problemformulation} formulates the
problem of interest and outlines the constraints underwhich our model is
applied. ~\autoref{sec:background} briefly reviews the core concepts of deep learning
used in this work, including recurrent and convolutional networks. A brief
review of projection-based model reduction is algo given in this section.
Finally, we review the important connection between autoencoders and POD. The
key contributions of this work are presented in ~\autoref{sec:nMOR}. Namely, the
construction of the convolutional autoencoder for nonlinear dimensionality
reduction and the construction of our modified LSTM network for modeling of
feature dynamics. In this section, we also discuss the construction of the
training datasets and develop our training strategy. ~\autoref{sec:numexp}
demonstrates the use of our method on three illustrative examples. The first
example considers a simple one-dimensional model reduction problem based on the
viscous Burgers equation. This serves to highlight the expressive power of
nonlinear autoencoders when compared to POD-based methods. Second, we consider a
parametric model reduction problem based on an incompressible flow inside a
periodic domain and evaluate our model's predictive performance with large
parameter variations. This example has the merit of showcasing the benefits of
the location-invariant properties of the convolutional autoencoder as compared
to POD-based models. The last example highlights the stability characteristics
of the convolutional recurrent autoencoder through a model reduction problem
based on a chaotic incompressible flow inside a lid-driven cavity. Finally,
~\autoref{sec:conclusions} presents a summary and discussion of our work.

\section{Problem forumaltion}
\label{sec:problemformulation}
\subsection{Nonlinear computational physics problem}
Consider a high-dimensional ODE resulting from the semi-discretization of a
time-dependent PDE
\begin{align}\label{eq:1}
  \begin{split}
    \dot{\mathbf{x}}(t) &= F(\mathbf{x}(t), t; \bm{\mu}),\\
    \mathbf{x}(t_0) &= \mathbf{x}_0(\bm{\mu}),
  \end{split}
\end{align}
where $t\in[t_0, T]\subset\mathbb{R}^+$ denotes time,
$\mathbf{x}\in\mathbb{R}^N$ is the spatially discretized state variable where
$N$ is large, and $\bm{\mu}\in\mathcal{D}\subseteq\mathbb{R}^d$ is the
vector of parameters sampled from the feasible parameter set $\mathcal{D}$.
Here, $F:\mathbb{R}^N\times\mathbb{R}^+\times\mathbb{R}^d\to\mathbb{R}^N$ is a
nonlinear function representing the dynamics of the discretized system. Such
large nonlinear systems are typical in the computational sciences such as when
numerically solving the Navier-Stokes equations describing a fluid flow. In the
parameter-varying case $\bm{\mu}$ may represent initial and boundary conditions,
material properties, or shape parameters of interest.

Often, in engineering design and analysis the interest is on the evolution of
certain outputs
\begin{equation}\label{eq:2}
  \mathbf{y} = G(\mathbf{x}(t),\bm{\mu}),
\end{equation}
where $\mathbf{y}\in\mathbb{R}^p$ may represent e.g., lift, drag, or some other
performance criteria. In this work, the attention is focused only on the
evolution of the full state $\mathbf{x}$.

\subsection{Completely data-driven model reduction}
When the number of degrees of freedom $N$ is large, evaluating ~\autoref{eq:1} for a
given initial condition and input parameter $\bm{\mu}$ becomes computationally
challenging in two particular applications. The first are \textit{time-critical}
applications, or applications where a solution needs to be attained within a
given threshold of time. Some examples include routine analysis applications and
model predictive control of distributed parameter systems where near-real time solutions are
crucial. The second are \textit{multi-query} applications, i.e., applications
where one needs to sample a large number of parameters from $\mathcal{D}$.
Examples of multi-query applications include shape optimization and uncertainty
quantification.

To alleviate this computational burden an offline-online strategy is usually
employed in which a dataset of solution snapshots
$\mathcal{X}=\{\mathbf{x}(t_i;\bm{\mu}_i)\}^{N_{\text{data}}}_{i=1}$ of
~\autoref{eq:1} is used to construct a surrogate model that is capable of
approximating new solutions at a fraction of the cost. A wide variety of
strategies exist for constructing these so-called data-driven models including
data-fit methods which use numerical or experimental data to fit an input-output
map, and projection-based reduced order models which approximately solve
~\autoref{eq:1} in a reduced subspace constructed from numerical or experimental
data.

While projection-based reduced order models are physics-based, and thus offer an
advantage over data-fit methods when it comes to physical interpretation, they
are often intrusive. That is one requires access to the operators when
performing the projection step. In this work, we will restrict our attention to
non-intrusive, purely data-driven reduced order modeling. Thus, the construction
of the surrogate model will require only the dataset $\mathcal{X}$ and no
information about ~\autoref{eq:1}. Indeed there are many situations, e.g. in
neuroscience and finance, in which data are abundant but governing laws are
uncertain or do not exist altogether. For the purposes of this work, we will
work under the assumption that we do not have access to ~\autoref{eq:1} from which
the datasets are generated.

\subsection{Single vs. multiple parameter-varying trajectories}
The construction and availability of solution snapshot datasets is inherently
problem dependent. Here, we focus on the two common cases encountered in model
reduction:
\begin{itemize}
  \item[(i)]The dataset $\mathcal{X} = \{\mathbf{x}(t_1;\bm{\mu}),\mathbf{x}(t_2;\bm{\mu}),...\}$
  is constructed using snapshots from a single, statistically stationary trajectory of
  ~\autoref{eq:1}. In this case, $\bm{\mu}$ is the same for all snapshots. This is
  relevant to situations in which obtaining snapshot data is exceedingly
  expensive such as in large direct numerical simulations and the interest is on
  obtaining ``quick'' approximate solutions.
  \item[(ii)]The dataset $\mathcal{X}=\{X^{\bm{\mu}_1}, X^{\bm{\mu}_{2}},...\}$
  is constructed using multiple, parameter varying trajectories
  $X^{\bm{\mu}_i}=\{\mathbf{x}(t_1;\bm{\mu}_i),\mathbf{x}(t_2;\bm{\mu}_i),...\}$. This case is
  relevant to multi-query applications or applications in which the interest
  is on capturing the parameter-dependent transient behavior of ~\autoref{eq:1}.
\end{itemize}
In both cases the surrogate model is constructed in a non-intrusive fashion
using the same procedure and only the dataset is changed.

\section{Background}
\label{sec:background}
In this section, we introduce the basic notions of deep learning and two key
architectures used in this work: 1) recurrent neural networks, and 2)
convolutional neural networks. Finally, we finish by summarizing the connections
between POD and fully-connected autoencoders.

\subsection{Deep learning}
Deep learning has enjoyed great success in recent years in areas from image and
speech recognition ~\cite{Krizhevsky2012,Couprie2013,Hinton2012} to genomics
~\cite{Xiong2015,Leung2014}. At the core of deep learning are deep neural
networks, whose layered structure allows them to learn at each layer a
representation of the raw input with increasing levels of abstraction
~\cite{Goodfellow2016a,LeCun2015a}. With enough layers, deep neural networks can
learn intricate structures in high-dimensional data. For example, given an image
as an array of pixel values, the first layer of a deep neural network might
learn to identify edges in various orientations. The second layer then is able
to detect particular arrangements of edges, and so on until a complex hierarchy
of features leads to the detection of a face or a road sign. Here, we briefly
review some concepts and common network architectures used in this work.

Neural networks are models of computation loosely inspired by biological
neurons. Generally, given a vector of real-valued inputs
$\mathbf{x}\in\mathbb{R}^N$, a single layer artificial neural network is an
affine transformation of the input $\mathbf{x}$ fed through a nonlinear function
\begin{equation}\label{eq:3}
  \hat{\mathbf{y}}=f(\mathbf{W}\mathbf{x}+\mathbf{b}),
\end{equation}
where $\mathbf{W}\in\mathbb{R}^{M\times N}$ is the weight matrix,
$\mathbf{b}\in\mathbb{R}^M$ is a bias term, and $f(\cdot)$
is a nonlinear function that acts element-wise on its inputs.

To create multilayered neural networks, the output $\mathbf{h}_l$ of a layer $l$
is fed as the input of the following layer, thus
\begin{align}\label{eq:4}
  \begin{split}
    \mathbf{h}_{l+1} &= f_{l+1}(\mathbf{W}_l\mathbf{h}_l+\mathbf{b}_l),\\
                     &= f_{l+1}(\mathbf{W}_lf_{l}(...(f_1(\mathbf{W}_0\mathbf{x}+\mathbf{b}_0))...+\mathbf{b}_{l-1})+\mathbf{b}_l),
  \end{split}
\end{align}
where $\mathbf{h}_1=f_1(\mathbf{W}_0\mathbf{x}+\mathbf{b}_0)$ is the output of
the first layer. The vector $\mathbf{h}_{l}$ is often referred to as the hidden
state or feature vector at the $l$-th layer. This process continues for $L$
layers, where at the final layer the output of the network is given by
$\hat{\mathbf{y}}=f_L(\mathbf{W}_{L-1}\mathbf{h}_{L-1}+\mathbf{b}_{L-1})$. In
supervised learning, training the network then involves finding the parameters
$\bm{\theta}=\{\mathbf{W}_l,\mathbf{b}_l\}_{l=0}^{L-1}$ such that the expected
loss between the output $\hat{\mathbf{y}}$ and the target value $\mathbf{y}$ is
minimized
\begin{equation}\label{eq:5}
  \bm{\theta}^*=\arg\min_{\bm{\theta}}\mathbb{E}_{(x,y)\sim\mathcal{P}_{data}}\big[ \mathcal{L}(f(\mathbf{x};\bm{\theta}),\mathbf{y}) \big],
\end{equation}
where $\mathcal{P}_{data}$ is the data-generating distribution and
$\mathcal{L}(\hat{\mathbf{y}},\mathbf{y})$ is some measure of discrepancy
between the predicted and target outputs. What distinguishes machine learning
from straight forward optimization is that the model $f$ parameterized by
$\bm{\theta}^*$ should be expected to generalize well for all examples drawn
from $\mathcal{P}_{data}$, even if they were not witnessed during training
~\cite{Goodfellow2016a}. Most neural networks are trained using stochastic
gradient descent (SGD), or one of its many variants
~\cite{Kingma2014,Zeiler2012}, in which gradients are computed using the
backpropagation procedure ~\cite{Rumelhart1986}.

\subsection{Recurrent neural networks}
A natural extension of feed-forward networks for sequential data are networks
with self-referential, or recurrent connection. These recurrent neural networks
(RNNs) process a sequence of inputs one element at a time, maintaining in the
hidden state an implicit history of previous inputs. Consider a sequence of
inputs $\{\mathbf{x}^{0}, \mathbf{x}^{1},...,\mathbf{x}^{m}\}$, with each
$\mathbf{x}^{n}\in\mathbb{R}^N$, the  $n$-th hidden state
$\mathbf{h}^n\in\mathbb{R}^{N_h}$ of a simple RNN is evaluated by the following
update
\begin{equation}\label{eq:6}
  \mathbf{h}^{n}=f(\mathbf{W}\mathbf{h}^{n-1}+\mathbf{U}\mathbf{x}^{n}+\mathbf{b}),
\end{equation}
where $\mathbf{W}\in\mathbb{R}^{N_h\times N_h}$ and
$\mathbf{U}\in\mathbb{R}^{N_h\times N}$ are the hidden and input weight matrices
respectively, and $\mathbf{b}\in\mathbb{R}^{N_h}$ is a bias term. RNNs are also
typically trained using SGD, or some variant, but the gradients are calculated
using the backpropagation through time (BPTT) algorithm ~\cite{Werbos1990}. In
BPTT, the RNN is first ``unrolled" in time, stacking one copy of the RNN per
time step. This results in a \textit{weight-tied} deep feed forward neural
network on which the standard backpropagation algorithm can be employed.

Training RNNs has long been considered to be challenging ~\cite{Bengio1994}. The
main difficulty is due to the exponential growth or decay of gradients as they
are backpropagated through each time step, so over many time steps they will
either vanish or explode. This is especially problematic when learning sequences
with long-term dependencies. The vanishing or exploding gradient problem is
typically addressed by using gated RNNs, including long short-term memory (LSTM)
networks ~\cite{Hochreiter1997} and networks based on the gated recurrent unit
(GRU) ~\cite{Cho2014a}. These networks have additional paths through which
gradients neither vanish nor explode, allowing gradients of the loss function to
backpropagate across multiple time-steps and thereby making the appropriate
parameter updates. This work will only consider RNNs equipped with LSTM units.

\subsection{Convolutional neural networks}
The final standard neural network architecture considered in this work are
convolutional neural networks. These networks were first introduced as an
alternative to fully connected networks for data structured as multiple arrays
(e.g., 1D signals and sequences, 2D images or spectrograms, and 3D video). The
two key properties of convolutional neural networks are: 1) local connections,
and 2) shared weights ~\cite{Goodfellow2016a,LeCun2015a}. In arrayed data, often
local groups of values are highly correlated, assembling into distinct features
that can be easily detected using a local approach. Additionally, weight sharing
across the input domain works to detect location-invariant features.

\begin{figure}[htbp]
  \centering
  \includegraphics[width=0.75\textwidth]{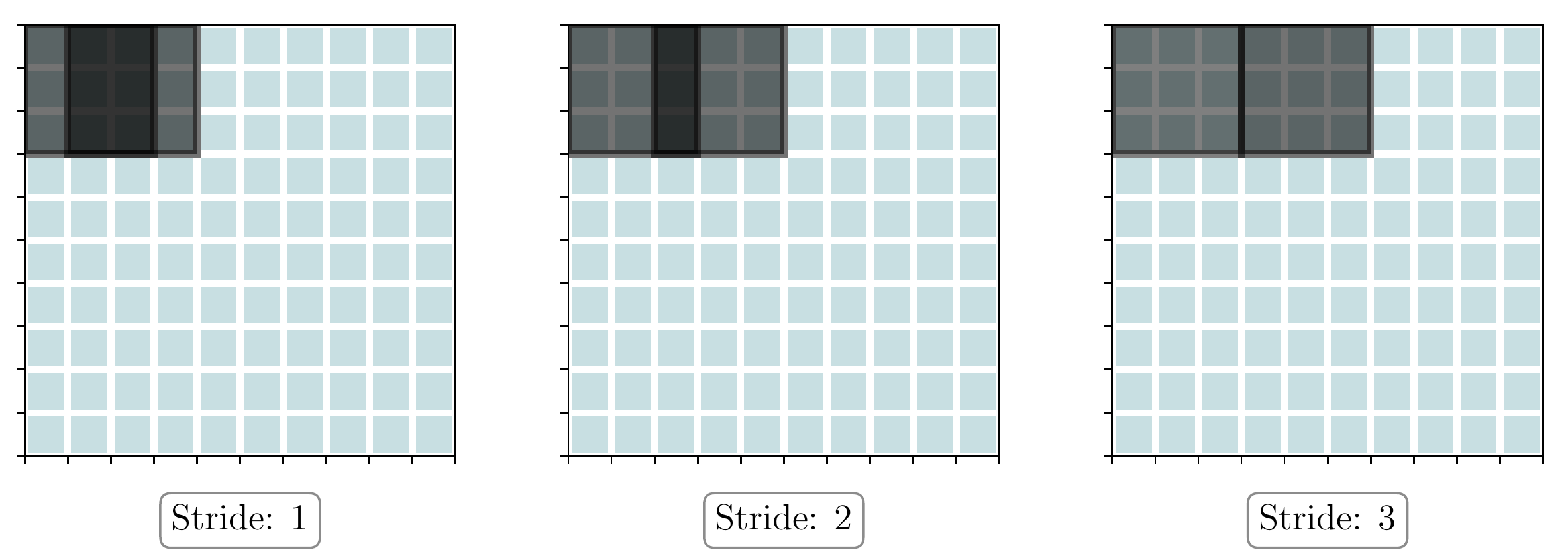}
  \caption{Sliding convolutional filter with varying stride values.}
  \label{fig:1}
\end{figure}

\begin{figure}[htbp]
  \centering
  \includegraphics[width=0.75\textwidth]{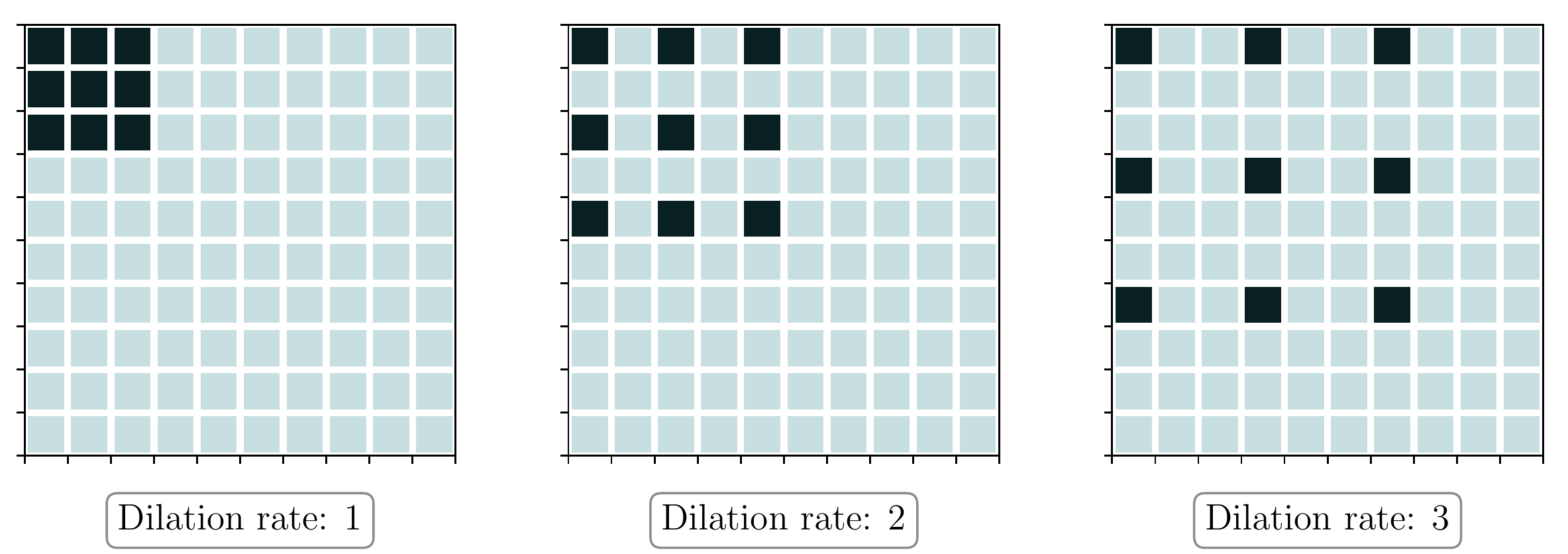}
  \caption{Convolutional filters with varying dilation rates.}
  \label{fig:2}
\end{figure}

In convolutional neural networks, layers are organized into feature maps, where
each unit in a feature map is connected to a local domain of the previous layer
through a filter bank. Consider a 2D input $\mathbf{X}\in\mathbb{R}^{N_x\times
N_y}$, a convolutional layer consists of a set of $F$ filters
$\mathbf{K}^f\in\mathbb{R}^{a\times b}$, $f=1,...,F$, each of which generates a
feature map $\mathbf{Y}^f\in\mathbb{R}^{N_x'\times N_y'}$ by a 2D discrete
convolution
\begin{equation}\label{eq:7}
  \mathbf{Y}^f_{i,j} = \sum^{a-1}_{k=0}\sum^{b-1}_{l=0}\mathbf{K}^f_{a-k,b-l}\mathbf{X}_{1+s(i-1)-k,1+s(j-1)-l},
\end{equation}
where $N_x'=1+\frac{N_x+a-2}{s}$, $N_y'=1+\frac{N_y+b-2}{s}$, and $s\ge1$ is an
integer value called the \textit{stride}. ~\autoref{fig:1} shows the effect of
different stride values of a filter acting on an input feature map. As before,
the feature map can be passed through an element-wise nonlinear function.
Typically, the dimension of the feature map is reduced by using a pooling layer,
in which a single value is computed from small $a'\times b'$ patch of the
feature map either by taking the maximum value or averaging. A slightly more
general approach is to employ a convolutional layer with a stride of $s>1$, in
which instead of taking the maximum or average value, some weighted sum of the
local patch of the input feature map is learned by adjusting the respective
filter $\mathbf{K}^f$. In addition, dilated convolutional filters (see
~\autoref{fig:2}) are often employed to significantly increase the receptive
field without loss of resolution, effectively capturing larger features in
highly dense data ~\cite{Yu2015, Li2018}.

\subsection{Projection-based model reduction}
In projection-based MOR, the state vector $\mathbf{x}\in\mathbb{R}^N$ is approximated
by a global affine trial subspace $\mathbf{x}^0+\mathcal{S}\subset\mathbb{R}^N$ of
dimension $N_h<<N$
\begin{equation}\label{eq:8}
  \mathbf{x}\approx\tilde{\mathbf{x}}=\mathbf{x}^0 + \bm{\Psi}_{N_h}\mathbf{h},
\end{equation}
where the columns of $\bm{\Psi}_{N_h}\in\mathbb{R}^{N\times N_h}$ contain the
basis for subspace $\mathcal{S}$, the initial condition is
given by $\mathbf{x}^0$, and $\mathbf{h}\in\mathbb{R}^{N_h}$ represents the
generalized coordinates in this subspace. Substituting ~\autoref{eq:8} into
~\autoref{eq:1} yields
\begin{equation}\label{eq:9}
  \bm{\Psi}_{N_h}\frac{d\mathbf{h}}{dt} = F(\mathbf{x}^0 + \bm{\Psi}_{N_h}\mathbf{h}(t; \bm{\mu})),
\end{equation}
which is an overdetermined system with $N$ equations and $N_h$ unknowns.
Additional constraints are imposed by enforcing the orthogonality of the
residual of ~\autoref{eq:9} on a test subspace represented by
$\bm{\mathbf{\Phi}}\in\mathbb{R}^{N\times N_h}$ through a Petrov-Galerkin
projection
\begin{equation}\label{eq:10}
  \bm{\Phi}^TR(\mathbf{x}^0 + \bm{\Psi}_{N_h}\mathbf{h}(t; \bm{\mu})) = 0,
\end{equation}
resulting in a square system with $N_h$ equations and $N_h$ unknowns, where
$R(\cdot)$ represents the residual of ~\autoref{eq:9}. In a Galerkin projection,
$\bm{\Phi}=\bm{\Psi}_{N_h}$. An important task is now choosing the subspace $\mathcal{S}$ on
which to project the governing equation ~\autoref{eq:1}. One popular method is to
obtain the basis of $\mathcal{S}$ through proper orthogonal decomposition
(POD).\footnotemark \footnotetext{This method is known under different names in
various fields: POD, principal component analysis (PCA), Karhunen–Lo\`eve
decomposition, empirical orthogonal functions and many others. In this work we
will adopt the name POD.}

Beginning with a set of $m$ observations $\{\mathbf{x}^n\}^m_{n=1}$,
$\mathbf{x}^n\in\mathbb{R}^N$, formed into a data matrix
$\mathbf{X}=[\mathbf{x}^1, \mathbf{x}^2, ...,
\mathbf{x}^m]\in\mathbb{R}^{N\times m}$, POD consists of performing the singular
value decomposition (SVD) on this data matrix
\begin{equation}\label{eq:11}
  \mathbf{X}=\mathbf{\Psi}\mathbf{\Sigma}\mathbf{V}^T,
\end{equation}
where $\mathbf{\Psi}\in\mathbb{R}^{N\times r}$ and
$\mathbf{V}\in\mathbb{R}^{N\times r}$ are orthonormal, i.e.,
$\mathbf{\Psi}^T\mathbf{\Psi}=\mathbf{V}^T\mathbf{V}=\mathbf{I}_{r\times r}$,
and $\mathbf{\Sigma}\in\mathbb{R}^{r\times r}$ is a diagonal matrix of whose
entries $\sigma_i\ge0$, ordered as $\sigma_1\ge\sigma_2\ge...\ge\sigma_r$ are
the singular values. The columns $\bm{\psi}_i$ of $\mathbf{\Psi}$ are sometimes called
the principal components, features, or POD modes. These modes have the property
that the linear subspace $\mathcal{S}$ spanned by
$\bm{\Psi}_{N_h}=[\bm{\psi}_1,...,\bm{\psi}_{N_h}]$, $N_h<r$, optimally
represents the data in the $L_2$ sense
\begin{equation}\label{eq:12}
  \min_{\mathbf{\Psi}_{N_h}}\| \mathbf{X} - \mathbf{\Psi}_{N_h}\mathbf{\Psi}^T_{N_h}\mathbf{X} \|^2_2 =
  \min_{\mathbf{\Psi}_{N_h}}\sum^{m}_{i=1}\| \mathbf{x}_i - \mathbf{\Psi}_{N_h}\mathbf{\Psi}^T_{N_h}\mathbf{x}_i \|^2_2.
\end{equation}
The net result is an optimal low-dimensional representation
$\mathbf{h}=\mathbf{\Psi}^T_{N_h}\mathbf{x}$ of an input $\mathbf{x}$, where
again $\mathbf{h}$ can be thought of as the intrinsic coordinates on the linear
subspace $\mathcal{S}$.

\subsection{Connection between autoencoders and POD}
In data-driven sciences, dimensionality reduction attempts to approximately
describe high-dimensional data in terms of a low-dimensional representation.
Central to this is the \textit{manifold hypothesis}, which presumes that
real-world high-dimensional data lies near a low-dimensional manifold
$\mathcal{S}$ embedded in $\mathbb{R}^N$, where $N$ is large ~\cite{Bengio2013}.
As a result POD has found broad applications from pre-training machine learning
models to dimensionality reduction of physical systems. However it has the major
drawback of constructing only an optimal linear manifold. This is quite
significant since data sampled from complex, real-world systems is more often
than not strongly nonlinear. A wide variety of strategies for more accurate
modeling of $\mathcal{S}$ have been developed over the years, most involving
using a patchwork of local subspaces $\{\mathcal{S}_l\}^L_{l=1}$ obtained
through linearizations or higher-order approximations of the state-space
~\cite{Rewienski2006,Bengio2013,Trehan2016}.

A nonlinear generalization of POD is the under-complete autoencoder
~\cite{Hinton2006,Goodfellow2016a}. An under-complete autoencoder consists of a
single or multiple-layer \textit{encoder} network
\begin{equation}\label{eq:13}
  \mathbf{h}=f_E(\mathbf{x};\bm{\theta}_E),
\end{equation}
where $\mathbf{x}\in\mathbb{R}^N$ is the input state,
$\mathbf{h}\in\mathbb{R}^{N_h}$ is the feature or representation vector, and
$N_h<N$. A \textit{decoder} network is then used to reconstruct $\mathbf{x}$ by
\begin{equation}\label{eq:14}
  \hat{\mathbf{x}} = f_D(\mathbf{h};\bm{\theta}_D).
\end{equation}
Training this autoencoder then consists of finding the parameters that minimize
the expected reconstruction error over all training examples
\begin{equation}\label{eq:15}
  \bm{\theta}^*_E,\bm{\theta}^*_D = \arg\min_{\bm{\theta}_E,\bm{\theta}_D} \mathbb{E}_{x\sim\mathcal{P}_{data}}\big[ \mathcal{L}(\hat{\mathbf{x}},\mathbf{x}) \big],
\end{equation}
where $ \mathcal{L}(\hat{\mathbf{x}},\mathbf{x})$ is some measure of discrepancy
between $\mathbf{x}$ and its reconstruction $\hat{\mathbf{x}}$. Restricting
$N_h<N$ serves as a form of regularization, preventing the autoencoder from
learning the identify function. Rather, it captures the salient features of the
data-generating distribution $\mathcal{P}_{data}$. Under-complete autoenocders
are just one of a family of regularized autoencoders which also include
contractive autoencoders, denoising autoencoders, and sparse autoencoders
~\cite{Hinton2006, Goodfellow2016a}.

\begin{remark}
  The choice of $f_E$, $f_D$, and $\mathcal{L}(\hat{\mathbf{x}},\mathbf{x})$
  largely depends on the application. Indeed, if one chooses a linear encoder
  and a linear decoder of the form
  \begin{align}\label{eq:16}
    \mathbf{h} &= \mathbf{W}_E\mathbf{x},\\
    \hat{\mathbf{x}} &= \mathbf{W}_D\mathbf{h},
  \end{align}
  where $\mathbf{W}_E\in\mathbb{R}^{N_h\times N}$ and
  $\mathbf{W}_D\in\mathbb{R}^{N\times N_h}$, then with a squared reconstruction error
  \begin{align}\label{eq:17}
  \begin{split}
    \mathcal{L}(\hat{\mathbf{x}}, \mathbf{x}) &= \| \mathbf{x} - \hat{\mathbf{x}}\|^2_2\\
                                              &= \| \mathbf{x} - \mathbf{W}\mathbf{W}^T\mathbf{x}\|^2_2,
  \end{split}
  \end{align}
  the autoencoder will learn the same subspace as the one spanned by the first $N_h$
  POD modes if
  $\mathbf{W}=\mathbf{W}_D=\mathbf{W}_E^T$. However, without additional
  constraints on $\mathbf{W}$, i.e.,
  $\mathbf{W}^T\mathbf{W}=\mathbf{I}_{N_h\times N_h}$, the columns of
  $\mathbf{W}$ will not form an orthonormal basis or have any hierarchical
  ordering ~\cite{Bengio2013,Plaut2018}.
\end{remark}

\section{Convolutional recurrent autoencoders for model reduction}
\label{sec:nMOR}

\subsection{Previous work and objectives}
In the past few years machine learning has become an increasingly attractive
tool in modeling or augmenting low-dimensional models of complex systems.
Broadly, machine learning has been used in three ways in this respect: 1) as
input-output maps to model closure terms in unstable POD-Galerkin models, 2) as
a means to model the evolution of the intrinsic coordinates from an optimal
subspace approximation of the state, and more recently 3) as an approach to
construct end-to-end models that both find optimal representations of the system
variables and linearly evolve these representations. Our work was motivated and
thus has important similarities to previous work in both the second and third
approach.

The main idea behind modeling the evolution of the optimal subspace
approximations of the state variable directly addresses one of the main
challenges of projection-based model reduction. Namely, for systems where
governing laws do not exist, a simple yet powerful approach is to model the
evolution the intrinsic coordinates, obtained for example through POD, using
of a recurrent neural network
\begin{equation}\label{eq:18}
  \mathbf{h}^{n+1} = f_{RNN}(\mathbf{h}^{n}),
\end{equation}
where the representation vector $\mathbf{h}\in\mathbb{R}^{N_h}$, is of much
lower dimension than the data from which it is approximated. This strategy has
previously been explored in the context of model reduction where $\mathbf{h}$ is
obtained through POD ~\cite{Kani2017a, Wang2017a}, and in the more general case
where ~\autoref{eq:18} may model the dynamic behavior of complex ~\cite{Ogunmolu2016}
or chaotic systems ~\cite{Yeo2017a}. This opens up a family of strategies for
modeling the dynamics of not just systems without HFM, but systems with
heterogeneous data sources, and systems with \textit{a priori} unknown optimal
subspace approximations -- a feature which we make use of in this work.

A more completely data-driven approach, and one that is more closely related
to our work, is to both learn a low-dimensional representation of the state
variable and to learn the evolution of this representation. This approach has
been explored in ~\cite{Otto2017a} in which an autoencoder is used to learn
a low-dimensional representation of the high-dimensional state,
\begin{equation}\label{eq:19}
  \mathbf{h} = f_E(\mathbf{x}),
\end{equation}
where $\mathbf{x}\in\mathbb{R}^N$ high-dimensional state of the system,
$\mathbf{h}\in\mathbb{R}^{N_h}$, $N_h<N$, and a linear recurrent model is used
to evolve the low-dimensional features
\begin{equation}\label{eq:20}
  \mathbf{h}^{n+1} = \mathbf{K}\mathbf{h}^{n},
\end{equation}
where $\mathbf{K}\in\mathbb{R}^{N_h\times N_h}$. This approach was first
introduced in the context of learning a dictionary of functions used in extended
dynamic mode decomposition to approximate the Koopman operator of a nonlinear
system ~\cite{Li2017b}.

The central theme in these approaches and projection-based model reduction in
general is the following two-step process:
\begin{enumerate}
  \item The identification of a low-dimensional manifold $\mathcal{S}$ embedded
  in $\mathbb{R}^N$ on which most of the data is supported. This yields, in
  some sense, an optimal low-dimensional representations $\mathbf{h} =
  f(\mathbf{x})$ of the data $\mathbf{x}$ in terms of intrinsic coordinates on
  $\mathcal{S}$, and
  \item The identification of a dynamic model which efficiently evolves the low-dimensional
  representation $\mathbf{h}$ on the manifold $\mathcal{S}$.
\end{enumerate}
In this work, we build on the framework introduced in ~\cite{Wang2017a,
Kani2017a, Otto2017a} for constructing or augmented reduced order models, and
extend it in multiple directions. First, we introduce a deep convolutional
autoencoder architecture which provides certain advantages in identifying
low-dimensional representation of the input data. Second, since the dynamics of
reduced state vector on $\mathcal{S}$ may not necessarily be linear, we employ a
single-layer LSTM network to model the possibly nonlinear evolution of
$\mathbf{h}$ on $\mathcal{S}$. Lastly, we introduce an unsupervised training
strategy which trains the convolutional autoencoder while using the current
reduced state vectors to dynamically train the LSTM network.

\subsection{Dimensionality reduction via convolutonal autoencoders}
Dimensionality reduction through fully-connected autoencoders have long been
used in a wide variety of applications. However, one quickly runs into the
\textit{curse of dimensionality} when considering DNS-level input data which can
easily reach $10^9$ dof as mentioned in the introduction. Directly applying
large physical or simulation data to fully-connected autoencoders is not only
computationally prohibitive, but the approach itself ignores the opportunity to
exploit the structure of features in high-dimensional data. That is, since
fully-connected autoencoders require that the input data be flattened into an 1D
array, the local spatial relations between values are eliminated and can only be
recovered by initially considering dense models. Sparsity can be achieved either
\textit{a posteriori} by pruning individual connections (setting $w_{ij}=0$ for
some $i,j$) after training, or encouraged during training by using $L_1$
regularization. Here, we seek to exploit local correlations present in many
physics-based data through the use of convolutional neural networks.

\begin{figure}[htbp]
  \centering
  \includegraphics[width=0.85\textwidth]{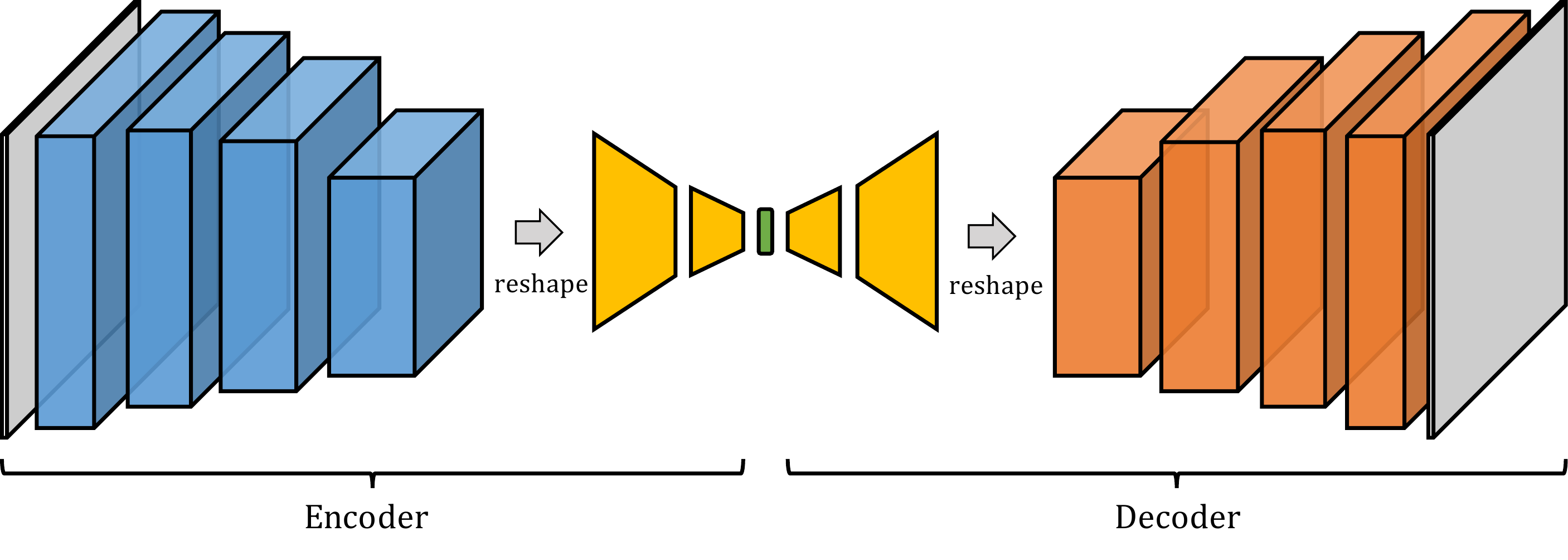}
  \caption{Network architecture of the convolutional autoencoder. The encoder
  network consists of a 4-layer convolutional encoder (blue), a 4-layer
  fully-connected encoder and decoder (yellow), and a 4-layer transpose
  convolutional decoder (red). The low-dimensional representation is depicted in
  green.}
  \label{fig:5}
\end{figure}

In particular, rather than applying a fully-connected autoencoder directly to
complex, high-dimensional simulation or experimental data, we apply it to a
vectorized feature map of a much lower-dimension obtained from a deep
convolutional network acting directly on the high-dimensional data. Wrapping the
fully-connected autoencoder with a convolutional neural network has two
significant advantages:
\begin{itemize}
  \item[(i)]{The local approach of each convolutional layer helps to exploit
  local correlations in field values. Thus, much the same way finite-difference
  stencils can capture local gradients, each filter $\mathbf{K}^f$ in a filter
  bank computes local low-level features from a small subset of the input.}
  \item[(ii)]{The shared nature of each filter bank both allows to identify
  similar features throughout the input domain and reduce the overall number of
  trainable parameters compared to a fully-connected layer with the same input
  size.}
\end{itemize}

Consider the following 12-layer convolutional autoencoder model depicted
graphically in ~\autoref{fig:5}. A 2D arrayed input
$\mathbf{X}\in\mathbb{R}^{N_x\times N_y}$, with $N_x=N_y=128$, is first passed
through 4-layer convolutional encoder. Each convolutional encoder layer uses a
filter bank $\mathbf{K}^f\in\mathbb{R}^{5\times 5}$, with the first layer having
a dilation rate of 2 and the number of filters $f$ increasing from 4 in the
first layer to 32 in the fourth layer using ~\autoref{eq:7}. At the opposite end of
the convolutional autoencoder network we use a 4-layer decoder network
consisting of \textit{transpose} convolutional layers. Often erroneuously
referred to as ``deconvolutional" layers, transpose convolutional layers
multiply each element of the input with a filter $\mathbf{K}^f$ and sum over the
resulting feature map, effectively swapping the forward and backward passes of a
regular convolutional layer. The effect of using transpose convolutional layers
with a stride $s>1$ is two decode low-dimensional abstract features to a larger
dimensional representation. ~\autoref{table:1} outlines the architecture of the
convolutional encoder and decoder subgraphs. In this work will consider the
sigmoid activation function $\sigma(s)=1/1+\exp(-s)$ for each layer of the
autoencoder.\footnotemark \footnotetext{In recent years the rectified linear
units (ReLUs) ~\cite{Goodfellow2016a}, given by $\text{ReLU}(s)=\max(0,s)$, and
its many variants like the ELUs ~\cite{Clevert2015}, have been favored over the
sigmoid activation function. However, in this work we have found that ReLUs
produce results similar to linear model reduction theory since $\text{ReLU}(s)$
are linear for inputs $s\in\mathbb{R}^+$.}

\begin{table}[h!]
\centering
\caption{Convolutional encoder (left) and decoder (right) filter sizes and strides}
\begin{tabular}{||c c c c||}
 \hline
 Layer & filter size & filters & stride \\ [0.5ex]
 \hline\hline
 1 & $5\times 5$ & 4 & $2\times 2$ \\
 2 & $5\times 5$ & 8 & $2\times 2$ \\
 3 & $5\times 5$ & 16 & $2\times 2$ \\
 4 & $5\times 5$ & 32 & $2\times 2$ \\
 \hline
\end{tabular}
\begin{tabular}{||c c c c||}
 \hline
 Layer & filter size & filters & stride \\ [0.5ex]
 \hline\hline
 9 & $5\times 5$ & 16 & $2\times 2$ \\
 10 & $5\times 5$ & 8 & $2\times 2$ \\
 11 & $5\times 5$ & 4 & $2\times 2$ \\
 12 & $5\times 5$ & 1 & $2\times 2$ \\
 \hline
\end{tabular}
\label{table:1}
\end{table}

In between the convolutional encoder and decoder is a regular fully-connected
autoencoder consisting of a 2-layer encoder which takes the vectorized form of
the $32$ feature maps from the last convolutional encoder layer
$\text{vec}(\mathcal{Y})\in\mathbb{R}^{512}$, where
$\mathcal{Y}=[\mathbf{Y}^1,...,\mathbf{Y}^{32}]\in\mathbb{R}^{4\times 4\times 32}$,
and returns the final low-dimensional representation of the input data
\begin{equation}\label{eq:21}
  \mathbf{h}=\sigma(\mathbf{W}^2_E\sigma(\mathbf{W}^1_E\text{vec}(\mathcal{Y})+\mathbf{b}^1_E)+\mathbf{b}^2_E)\in\mathbb{R}^{N_h},\quad N_h<<(N_x\cdot N_y),
\end{equation}
where $\mathbf{W}^1_E,\mathbf{b}^1_E$ and $\mathbf{W}^2_E,\mathbf{b}^2_E$ are
the parameters of the first and second fully-connected encoder network (the
$5^\text{th}$ and $6^\text{th}$ layers of the whole model). To reconstruct the
original input data from the low-dimensional representation, a similar 2-layer
fully-connected decoder parameterized by $\mathbf{W}^1_D,\mathbf{b}^1_D$ and
$\mathbf{W}^2_D,\mathbf{b}^2_D$, whose result is reshaped and passed to the
transpose convolutional decoder network.

Hierarchical convolutional feature learning through similar strategies have
previously been proposed for visual tracking ~\cite{Ma2015} and scene labeling or
semantic segmentation of images ~\cite{Couprie2013,Noh2015,Badrinarayanan2017}.
However, this is the first time, to the authors knowledge, that a convolutional
autoencoders have been applied to model reduction of large numerical data of
physical dynamical systems. The key innovation in using convolutional
autoencoders in model reduction is that it allows for nonlinear autoencoders and
thus nonlinear model reduction to be applied to large input data in a way that
exploits structures inherent in many physical systems.

\begin{remark}
  In this work we restrict our attention to 2D input data of size $N_x\times
  N_y=128\times 128$ with the first layer convolutional filter having a dilation
  rate of 2. In practice, however, an equivalent memory-reducing approach was
  employed by using an input data of size $N_x\times N_y=64\times 64$. In
  addition, the low-dimensional representations considered in this work are of
  size $N_h=64$ or smaller. To this effect, the hidden state sizes of the middle
  fully-connected autoencoder were chosen to be $512$ and $256$ such that
  $\mathbf{W}^1_E,(\mathbf{W}^2_D)^T\in\mathbb{R}^{512\times 256}$ and
  $\mathbf{W}^2_E,(\mathbf{W}^1_D)^T\in\mathbb{R}^{256\times N_h}$ with the bias
  terms shaped accordingly. The net result is an autoencoder with a maximum of
  330k parameters with $N_h=64$. A similar 12-layer fully-connected autoencoder
  would require over 22M parameters.
\end{remark}

\begin{remark}
  The size of the low-dimensional representation $N_h$ must be chosen
  \textit{a priori} for each model. Currently, no principled approach exists for
  the choice of $N_h$. One possible heuristic for an upper bound is to choose
  $N_h$ such that
  \begin{equation}
    \frac{\sum^{N_h}_{i=1}\sigma^2_i}{\sum^m_{i=1}\sigma^2_i}<\kappa,
  \end{equation}
  where $\sigma_i\ge0$ are the singular values of the data matrix
  $\mathbf{X}\in\mathbb{R}^{N\times m}$ and $\kappa$ is usually taken to be
  $99.9\%$. This approach is often employed when selecting the number of POD
  modes to keep in POD-Galerkin reduced order models ~\cite{Benner2015}, where in
  the context of fluid flows this corresponds to choosing enough modes such that
  $99.9\%$ of the energy content in the flow is preserved.
\end{remark}

\subsection{Learning feature dynamics}
The second component of projection-based model reduction is modeling the
evolution of low-dimensional features $\mathbf{h}$ in a computationally
efficient manner. Though identifying linear dynamics of $\mathbf{h}$ is
beneficial from an analysis perspective, here we will consider the general case
of learning arbitrary feature dynamics.

Consider a set of observations $\{\mathbf{x}^n\}^{m}_{n=0}$,
$\mathbf{x}^n\in\mathbb{R}^N$ obtained from a HFM or through experimental
sampling. Furthermore, for each observations consider some optimal
low-dimensional representation $\mathbf{h}^n\in\mathbb{R}^{N_h}$, where
$N_h<<N$. This low-dimensional representation can come from an optimal
rank-$N_h$ POD representation $\mathbf{h}^n=\bm{\Psi}^T_{N_h}\mathbf{x}^n$,
where the columns of $\bm{\Psi}_{N_h}$ are the first $N_h$ POD modes, or through
a neural network approach such as an autoencoder. We seek to construct a model
for the evolution of this low-dimensional representation in a completely
data-driven fashion, i.e., without access or knowledge of the system operators.
This is particularly useful for cases where HFM are uncertain or do not exist
altogether.

To model the evolution of $\mathbf{h}$ we employ a modified version of the long
short term memory (LSTM) network. LSTM networks we first proposed primarily to
overcome the vanishing or exploding gradient problem and are equipped with an
explicit memory cell and four gating units which adaptively control the flow of
information through the network ~\cite{Hochreiter1997, Lipton2015}. LSTM networks
have demonstrated impressive results in modeling relationships between sequences
such as in machine translation tasks ~\cite{Sutskever2014a, Cho2014a, Wu2016,
LeCun2015a}. More recently, they have been used to predict conditional
probability distributions in chaotic dynamical systems ~\cite{Yeo2017a} and in
modeling the evolution of low-dimensional POD representations
~\cite{Wang2017a, Kani2017a}.

In this work we are interested in evolving feature vectors whose size
correspond to the intrinsic dimensionality of a physical system, which may be
small compared to the number of hidden states and layers used in e.g., machine
translation. In addition, for large-scale systems it may be inefficient, if not
computationally prohibitive to reconstruct the full high-dimensional state at
every time-step. With these restrictions, we construct a modified single-layer
LSTM network to evolve the low-dimensional representation $\mathbf{h}^n$ without
full state reconstruction with the following components:
\begin{itemize}
  \item{\textit{Input gate}:
  \begin{equation*}
    \mathbf{i}^n = \sigma(\mathbf{W}_i\mathbf{h}^{n-1}+\mathbf{b}_i)
  \end{equation*}}
  \item{\textit{Forget gate}:
  \begin{equation*}
    \mathbf{f}^n = \sigma(\mathbf{W}_f\mathbf{h}^{n-1}+\mathbf{b}_f)
  \end{equation*}}
  \item{\textit{Output gate}:
  \begin{equation*}
    \mathbf{o}^n = \sigma(\mathbf{W}_o\mathbf{h}^{n-1}+\mathbf{b}_o)
  \end{equation*}}
  \item{\textit{Cell state}:
  \begin{equation*}
    \mathbf{c}^n = \mathbf{i}^n\odot\mathbf{c}^{n-1} + \mathbf{i}^n\odot\tanh(\mathbf{W}_c\mathbf{h}^{n-1}+\mathbf{b}_c)
  \end{equation*}}
\end{itemize}
where all four gates are used to update the feature vector by
\begin{equation}\label{eq:22}
  \mathbf{h}^n = \mathbf{o}^n\odot\tanh(\mathbf{c}^n).
\end{equation}
Here, $\odot$ represents the Hadamard product. Intuitively, at each step $n$ the
input and forget gates choose what information gets passed and dropped from the
cell state $\mathbf{c}^n$, while the output gate controls the flow of
information from the cell state to the feature vector. It is important to note
that the the evolution of $\mathbf{h}$ does not require information from the
full state $\mathbf{x}$, thereby avoiding a costly reconstruction at every step.

\begin{figure}
  \centering
  \includegraphics[width=0.6\textwidth]{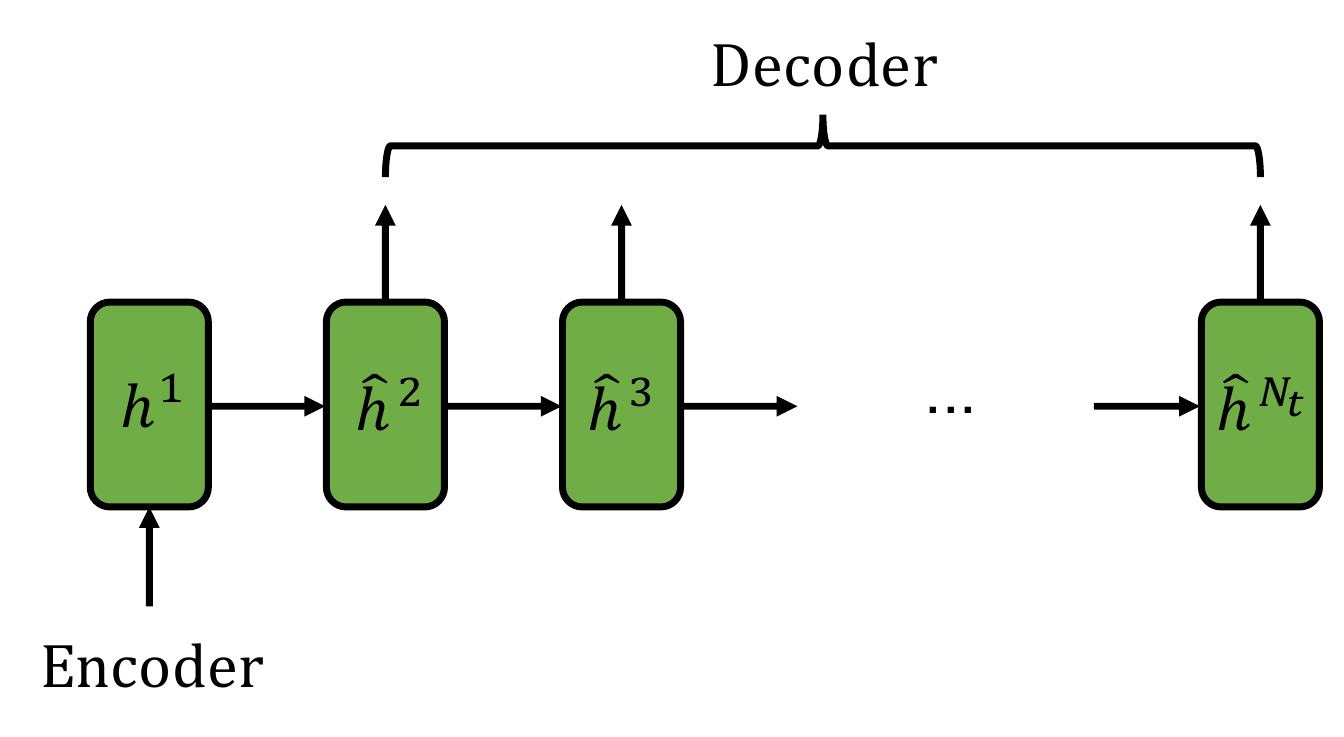}
  \caption{LSTM model that iteratively updates the low-dimensional representation $\mathbf{h}$.}
  \label{fig:6}
\end{figure}

Initializing with a known low-dimensional representation $\mathbf{h}^0$ one
obtains a prediction for the following steps by iteratively applying
\begin{equation}\label{eq:23}
  \mathbf{\hat{h}}^{n+1} = f_{LSTM}(\mathbf{\hat{h}}^n)\quad n=1,2,3,...
\end{equation}
where $\mathbf{\hat{h}}^1 = f_{LSTM}(\mathbf{h}^0)$, and $f_{LSTM}(\cdot)$
represents the action of ~\autoref{eq:22} and its subcompoents. A graphical
representation of this model is depicted in ~\autoref{fig:6}.

\subsection{Unsupervised training strategy}
A critical component of this work is the development of an unsupervised training
approach that adjusts both the convolutional autoencoder and recurrent model in a
joint fashion. The main obstacle is in preventing either the convolutional
autoencoder or RNN portion of the model from overfitting. Here, we discuss the
construction of the training dataset as well as the training and evaluation
algorithms.

\subsubsection{Constructing the training dataset}
Consider a dataset $\{\mathbf{x}^1,\mathbf{x}^2,...,\mathbf{x}^m\}$, where
$\mathbf{x}\in\mathbb{R}^{N_x\times N_y}$ is a 2D snapshot of some dynamical
system (e.g., a velocity field defined on a 2D grid). To make make this dataset
amenable to training it is broken up into a set of $N_s$ finite-time
training sequences $\{\mathbf{X}^1,...,\mathbf{X}^{N_s}\}$, where each training
sequence $\mathbf{X}^i\in\mathbb{R}^{N_x\times N_y\times N_t}$ consists of $N_t$
snapshots. Parameter-varying datasets naturally break up in this form where each
$\mathbf{X}^i$ may represent a small sequence of snapshots corresponding to a single
parameter value $\bm{\mu}_i$.

A common strategy for improving training is to consider only the fluctuations
around the temporal mean
\begin{equation}\label{eq:24}
  \mathbf{x}'^n = \mathbf{x}^n-\mathbf{\bar{x}},
\end{equation}
where $\mathbf{\bar{x}}=\frac{1}{m}\sum^m_{n=1}\mathbf{x}^n$ is the temporal
average over the entire dataset and $\mathbf{x}'$ are the fluctuations around
this mean. In our case, each layer in the convolutional autoencoder uses the
sigmoid activation function which maps each real-valued input to the interval
$(0,1)$, requiring our dataset be feature-scaled in order to prevent saturation
of the activation ~\cite{Goodfellow2016a}. Thus, our training dataset consists of
feature scaled snapshots
\begin{equation}\label{eq:25}
  \mathbf{x}'^n_s = \frac{\mathbf{x}'^n-\mathbf{x}'_{min}}{\mathbf{x}'_{max}-\mathbf{x}'_{min}},
\end{equation}
where each $\mathbf{x}'^n_s\in[0,1]^{N_x\times N_y}$. With these modifications,
the resulting training dataset has the following form
\begin{equation}\label{eq:26}
  \mathcal{X}=\{\mathbf{X}'^1_s,...,\mathbf{X}'^{N_s}_s\}\in[0,1]^{N_x\times N_y\times N_t\times N_s},
\end{equation}
where each training sample
$\mathbf{X}'^i_s=[\mathbf{x}'^1_{s,i},...,\mathbf{x}'^{N_t}_{s,i}]$ is a matrix
consisting of the feature-scaled fluctuations.

\subsubsection{Offline training and online prediction algorithms}
Our approach to train both components of the convolutional recurrent autoencoder
model is to split the forward pass into two stages. In the first stage, the
autoencoder takes an $N_b$-sized batch of the training data
$\mathcal{X}^b\subset\mathcal{X}$, where $\mathcal{X}^b\in[0,1]^{N_x\times
N_y\times N_t \times N_b}$, and outputs both the current $N_b$-sized batch of
low-dimensional representations of the training sequence
\begin{equation}\label{eq:27}
  \mathcal{H}^b=\{\mathbf{H}^1,...,\mathbf{H}^{N_b}\}\in\mathbb{R}^{N_h\times N_t\times N_b},
\end{equation}
where $\mathbf{H}^i=[\mathbf{h}^{1}_i,...,\mathbf{h}^{N_t}_i]\in\mathbb{R}^{N_h\times N_t}$ and a
reconstruction $\mathcal{\hat{X}}^b$ of the original input training batch. In
the second stage of the forward pass, the first feature vector of each sequence
is used to initialize and iteratively update ~\autoref{eq:23} to get a reconstruction
$\mathcal{\hat{H}}^b$ of the low-dimensional representations of the training batch ~\autoref{eq:26}.

We seek to construct a loss function that equally weights the error in the
full-state reconstruction and the evolution of the low-dimensional
representations. In general, we would like to find the model parameters
$\bm{\theta}$ such that for any sequence
$\mathbf{X}'_s=[\mathbf{x}'^1_{s},...,\mathbf{x}'^{N_t}_{s}]$, $\mathbf{x}'^n_{s}\sim\mathcal{P}_{data}$
and its corresponding low-dimensional representation
$\mathbf{H}=[\mathbf{h}^{1},...,\mathbf{h}^{N_t}]$, where $\mathcal{P}_{data}$
is the data-generating distribution, minimizes the following expected error between the model and the data
\begin{align}\label{eq:28}
  \begin{split}
    \mathcal{J}(\bm{\theta})  &= \mathbb{E}_{\mathbf{x}'^n_{s}\sim\mathcal{P}_{data}}
    \big[\mathcal{L}(\hat{\mathbf{X}}'_s,\mathbf{X}'_s,\hat{\mathbf{H}},\mathbf{H})\big]\\
                              &= \mathbb{E}_{\mathbf{x}'^n_{s}\sim\mathcal{P}_{data}}
    \bigg[
      \frac{\alpha}{N_t}\sum^{N_t}_{n=1}\frac{\|\mathbf{x}'^n_s-\hat{\mathbf{x}}'^n_s\|^2_F}{\|\mathbf{x}'^n_s\|^2_F+\epsilon} +
      \frac{\beta}{N_t-1}\sum^{N_t}_{n=2}\frac{\|\mathbf{h}^n-\hat{\mathbf{h}}^n\|^2_2}{\|\mathbf{h}^n\|^2_2+\epsilon}
    \bigg]
  \end{split}
\end{align}
where $\epsilon>0$ is a small positive number and $\alpha=\beta=0.5$. In
practice, the expected error is approximated by averaging
$\mathcal{L}(\hat{\mathbf{X}}'_s,\mathbf{X}'_s,\hat{\mathbf{H}},\mathbf{H})$
over all samples in a training batch during each backward pass. Intuitively, at
every training step, the autoencoder performs a regular forward pass while
constructing a new batch of low-dimensional representations which are used to
train the RNN. In this work we use the ADAM optimizer ~\cite{Kingma2014}, a
version of stochastic gradient descent that computes adaptive learning rates for
different parameters using estimates of first and second moments of the
gradients. Algorithm ~\ref{alg:1} outlines the  offline training of the
convolutional recurrent autoencoder in more detail. This model was built and
trained using the open-source deep learning library TensorFlow ~\cite{Abadi2016}.

\begin{algorithm2e}
\SetAlgoLined
\KwIn{Training dataset $\mathcal{X}\in[0,1]^{N_x\times N_y\times N_t\times N_s}$,
number of train-steps $N_{train}$, batch size $N_b$.}
\KwResult{Trained model parameters $\bm{\theta}$}
Randomly initialize $\bm{\theta}$\;
\For{$i\in\{1,...,N_{train}$\}}{
  Randomly sample batch from training data: $\mathcal{X}^{b}\subset\mathcal{X}$\;
  Flatten batch-mode: $\mathcal{X}^{b_{AE}}\leftarrow\text{flatten}(\mathcal{X}^{b})$ s.t. $\mathcal{X}^{b_{AE}}\in[0,1]^{N_x\times N_y\times (N_t\cdot N_b)}$\;
  Encoder forward pass: $\mathcal{\tilde{H}}^b\leftarrow f_{enc}(\mathcal{X}^{b_{AE}})$ where $\mathcal{\tilde{H}}^b\in\mathbb{R}^{N_h\times (N_t\cdot N_b)}$\;
  Decoder forward pass: $\hat{\mathcal{X}}^{b_{AE}}\leftarrow f_{dec}(\mathcal{\tilde{H}}^b)$\;
  Reshape low-dimensional features: $\mathcal{H}^b\in\mathbb{R}^{N_h\times N_t\times N_b}\leftarrow \text{reshape}(\mathcal{\tilde{H}}^b)$\;
  Initialize RNN subgraph loop: $\hat{\mathbf{h}}^2_i\leftarrow f_{LSTM}(\mathbf{h}^1_i)$ for $i\in\{1,...,N_b\}$, $\mathbf{h}_i^1\subset\mathcal{H}^b$\;
  \For{$n\in\{2,...,N_t-1\}$}{
    $\hat{\mathbf{h}}^{n+1}_i\leftarrow f_{LSTM}(\hat{\mathbf{h}}^n_i)$ for $i\in\{1,...,N_b\}$, $\hat{\mathbf{h}}^{n}_i\subset{\mathcal{\hat{H}}^b}$\;
  }
  Using $\mathcal{X}^b, \mathcal{\hat{X}}^b, \mathcal{H}^b$, and $\mathcal{\hat{H}}^b$ calculate approximate gradient $\hat{\mathbf{g}}$ of ~\autoref{eq:28}\;
  Update parameters: $\bm{\theta}\leftarrow ADAM(\hat{\mathbf{g}})$
}
\caption{Convolutional Recurrent Autoencoder Training Algorithm}
\label{alg:1}
\end{algorithm2e}

Once the model is trained, online prediction is straightforward. Using the
trained parameters $\bm{\theta}^*$, and given an initial condition
$\mathbf{x}^0\in[0,1]^{N_x\times N_y}$, a low-dimensional representation of the
initial condition $\mathbf{h}^0\in\mathbb{R}^N_h$ is constructed using the
encoder network. Iterative applications of ~\autoref{eq:23} are then used to evolve
this low-dimensional representation for $N_t$ steps. The modular construction of
the convolutional recurrent autoencoder model allows the user to reconstruct
from $\hat{\mathbf{h}}^n$ the full-dimensional state $\hat{\mathbf{x}}^n$ at
every time step or at any specific instance. The online prediction algorithm is
outlined in Algorithm ~\ref{alg:2}.

\begin{algorithm2e}
\SetAlgoLined
\KwIn{Initial condition $\mathbf{x}^0\in[0,1]^{N_x\times N_y}$, number of prediction steps $N_{t}$.}
\KwResult{Model prediction $\hat{\mathbf{X}}=[\hat{\mathbf{x}}^1,...,\hat{\mathbf{x}}^{N_t}]\in[0,1]^{N_x\times N_y\times N_t}$}
Load trained parameters $\bm{\theta}^*$\;
Encoder forward pass: $\mathbf{h}^0\leftarrow f_{enc}(\mathbf{x}^0)$\;
Initialize RNN subgraph loop: $\hat{\mathbf{h}}^1\leftarrow f_{LSTM}(\mathbf{h}^0)$\;
\For{$n\in\{1,...,N_t-1\}$}{
  $\hat{\mathbf{h}}^{n+1}\leftarrow f_{LSTM}(\hat{\mathbf{h}}^{n})$\;
}
Decoder forward pass: $\hat{\mathbf{X}}\leftarrow f_{dec}(\hat{\mathbf{H}})$, where $\hat{\mathbf{H}}=[\hat{\mathbf{h}}^1,...,\hat{\mathbf{h}}^{N_t}]$\;
\caption{Convolutional Recurrent Autoencoder Prediction Algorithm}
\label{alg:2}
\end{algorithm2e}

\section{Numerical experiments}
\label{sec:numexp}
We apply the methods described in the previous sections on three representative
examples to illustrate the effectiveness of deep autoencoder-based approaches to
nonlinear model reduction. The first one considers only a 4-layer
fully-connected recurrent autoencoder model applied to a simple one-dimensional
problem based on the viscous Burgers equation. This has the merit of
demonstrating the performance of autoencoders equipped with nonlinear activation
functions on tasks where linear methods tend to struggle. The second example
considers a parametric model reduction problem based on two-dimensional fluid
flow in a periodic domain with significant parameter variations. In this case,
our convolutional recurrent autoencoder model is tasked with predicting new
solutions given new parameters (i.e., parameters unseen during training). The
third example focuses on long-term prediction of an incompressible flow inside a
lid-driven cavity. This case serves to highlight the long-term stability and
overall performance of the convolutional recurrent autoencoder model in contrast
to the unstable behavior exhibited by POD-Galerkin ROMs.

\subsection{Viscous Burgers equation}
First, we consider the one-dimensional viscous Burgers equation given by
\begin{align}\label{eq:29}
\begin{split}
  &\frac{\partial u}{\partial t} + u\frac{\partial u}{\partial x} = \frac{1}{\text{Re}}\frac{\partial^2 u}{\partial x^2},\quad(x,t)\in[0,L]\times[0,T],\\
  &u(x,0) = 1 + \exp\bigg(-\frac{2(x-x_0)^2}{0.1^2}\bigg),\\
  &u(0,t) = 0,
\end{split}
\end{align}
where  $L=1.5$, $T=0.3$, $x_0$ is the initial location of the Gaussian initial condition,
and the Reynolds-like number is set to $Re=200$. This problem is
spatially discretized onto a uniform $N_x=1024$ grid using a second order finite
difference scheme with a grid spacing of $\Delta x=L/N_x$. A parameter-varying
dataset consisting of $N_s=128$ training samples is created by randomly sampling
$x_0\in[0,L]$ and solving ~\autoref{eq:29} using a fourth-order Runge-Kutta scheme
with $\Delta t=0.5\Delta x$. After subtracting the mean and feature scaling each
solution snapshot, the training dataset has the form of ~\autoref{eq:26} where each
training sample is a matrix of solution snapshots
$\mathbf{X}'^i_s=[\mathbf{u}'^1_{s,i},...,\mathbf{u}'^{N_t}_{s,i}]\in\mathbb{R}^{N_x\times
N_t}$ and corresponds to a different initial condition $x^i_0$. In this case,
$N_t=40$ is the number of equally spaced snapshots sampled from each trajectory.

This example was crafted to highlight an important performance benefit of using
nonlinear fully-connected autoencoder-based model reduction approaches in
contrast to POD-based ROMs. First, we train a 4-layer fully-connected
autoencoder (2 encoder layers and 2 decoder layers) to produce a low-dimensional representation
$\mathbf{h}_{AE}\in\mathbb{R}^{N_h}$, $N_h=20$ with intermediate layer of size
512. The evolution of this representation and an equivalently sized optimal POD
representation $\mathbf{h}_{POD}=\mathbf{\Psi}^T_{N_h}\mathbf{u}'_s$ are both
modeled using separate single layer modified LSTM networks trained according to
a simplified version of Algorithm ~\ref{alg:1} using a batch size $N_b=8$. A
best-case scenario for any POD-based ROM are snapshot reconstructions satisfying
~\autoref{eq:12}, therefore in lieu of a POD-Galerkin-ROM we will consider only the
projected solution snapshots. These proof-of-concept models were each trained
over $N_{train}=100,000$ iterations on a desktop computer in a matter of
minutes.

\begin{figure}[htbp]
  \centering
  \includegraphics[width=0.85\textwidth]{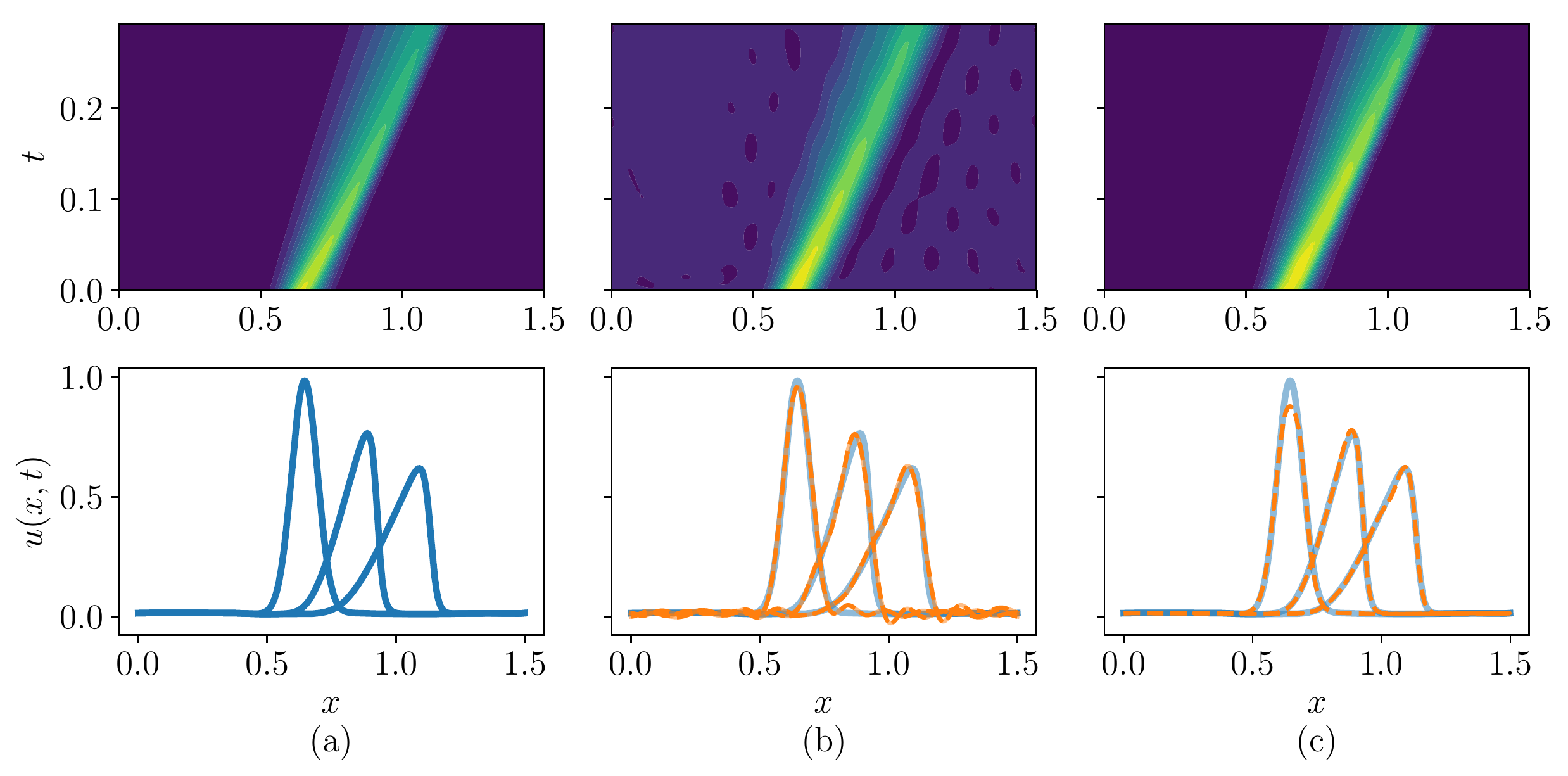}
  \caption{(a) exact solution, (b) $N_h=20$; $L_2$-optimal POD reconstruction
  (solid orange), POD-LSTM reconstruction (dashed orange), exact solution
  (light blue), (c) $N_h=20$; shallow autoencoder-LSTM reconstruction (dashed
  orange), exact solution (light blue).}
  \label{fig:7}
\end{figure}

~\autoref{fig:7} depicts the comparison between exact solution, the best-case optimal
POD reconstruction, the POD-LSTM reconstruction, and finally the shallow
recurrent autoencoder reconstruction. As expected, due to the truncation of
higher-frequency POD modes the $L_2$-optimal POD reconstruction exhibits
spurious oscillations. The spurious oscillations, aside from signifying a poor
reconstruction, may lead to stability issues in POD-Galerkin ROMs. This is
widely known to be a problem for model reduction of fluid flows where
POD-Galerkin ROMs, while capturing nearly all the energy of the system truncate
low-energy modes which can have a large influence on the dynamics. Additionally,
in agreement with similar work in ~\cite{Wang2017a, Kani2017a}, the POD-LSTM
model was able to accurately capture the evolution of the optimal POD
representation in a non-intrusive manner.

More importantly, the power of recurrent autoencoder-based approaches for
nonlinear model reduction is exhibited in the reconstruction using the shallow
recurrent autoencoder model. The effect of nonlinearities in the fully-connected
autoencoder help to identify a more expressive low-dimensional representation of
the full state. Combining this with an LSTM network to evolve these low-dimensional
representations yields an effective nonlinear reduced order modeling approach that
outperforms best-case scenario POD-based ROMs while using the same size models.

\subsection{Parameter-varying flow in a periodic box}
\begin{figure}[htbp]
  \centering
  \includegraphics[width=0.4\textwidth]{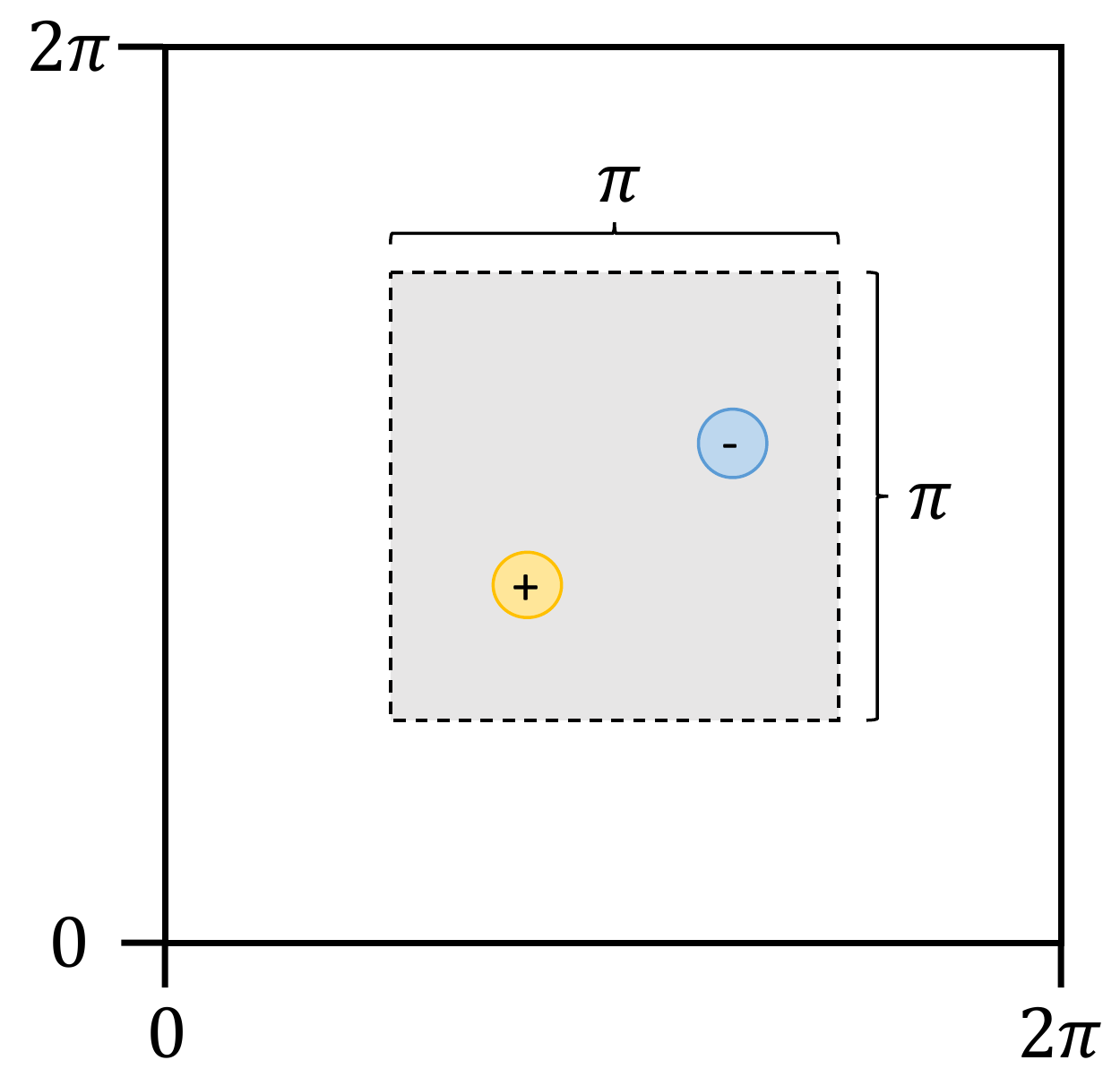}
  \caption{Square domain with periodic boundary conditions. The positive and negative
  vortices of equal strength are randomly initialized within the grey subdomain.}
  \label{fig:8}
\end{figure}

\begin{figure}[htbp]
  \centering
  \includegraphics[width=0.75\textwidth]{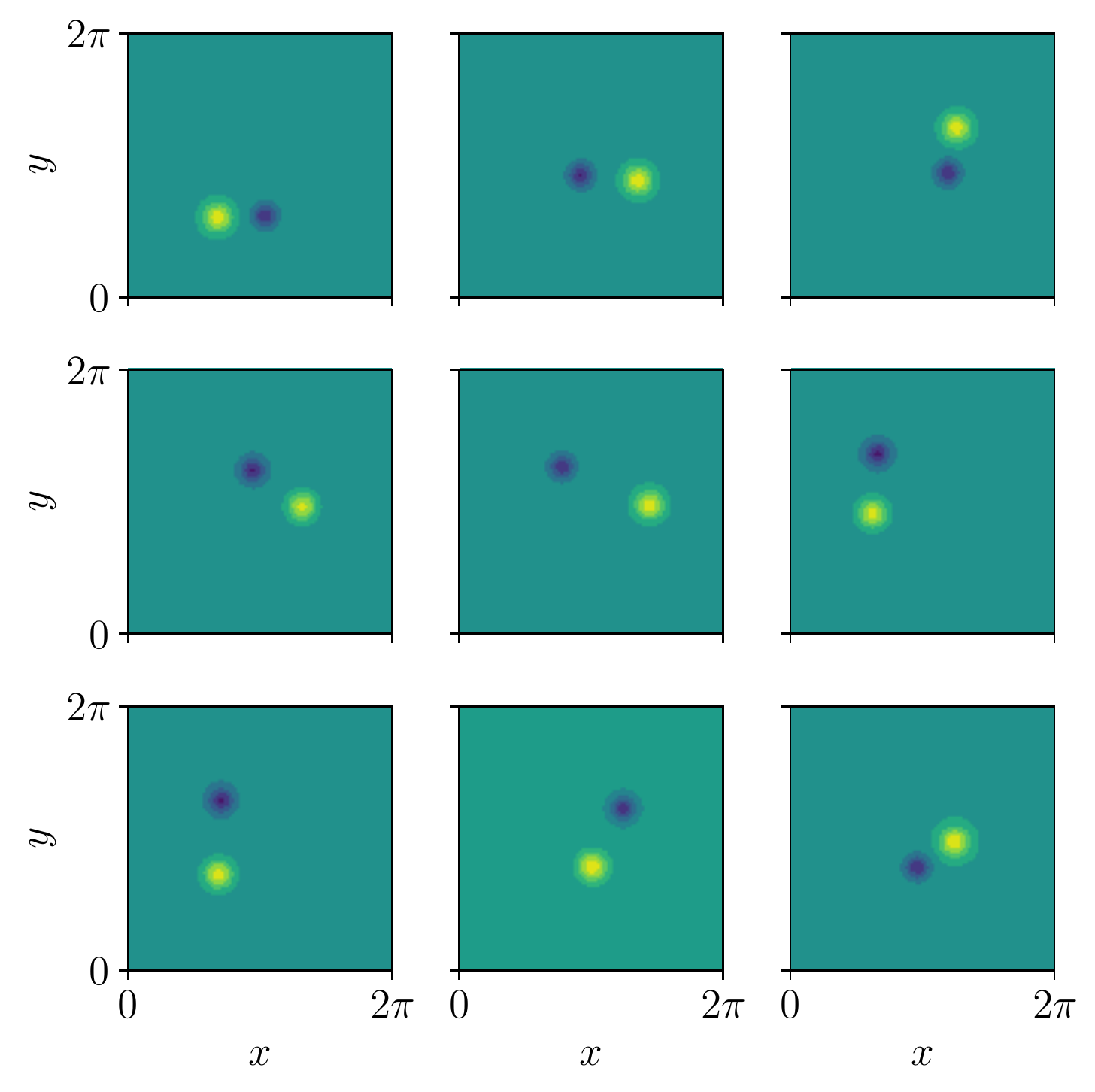}
  \caption{A set of initial conditions with two randomly located Gaussian vortices of
  equal and opposite strength.}
  \label{fig:9}
\end{figure}

\begin{figure}[htbp]
  \centering
  \includegraphics[width=0.75\textwidth]{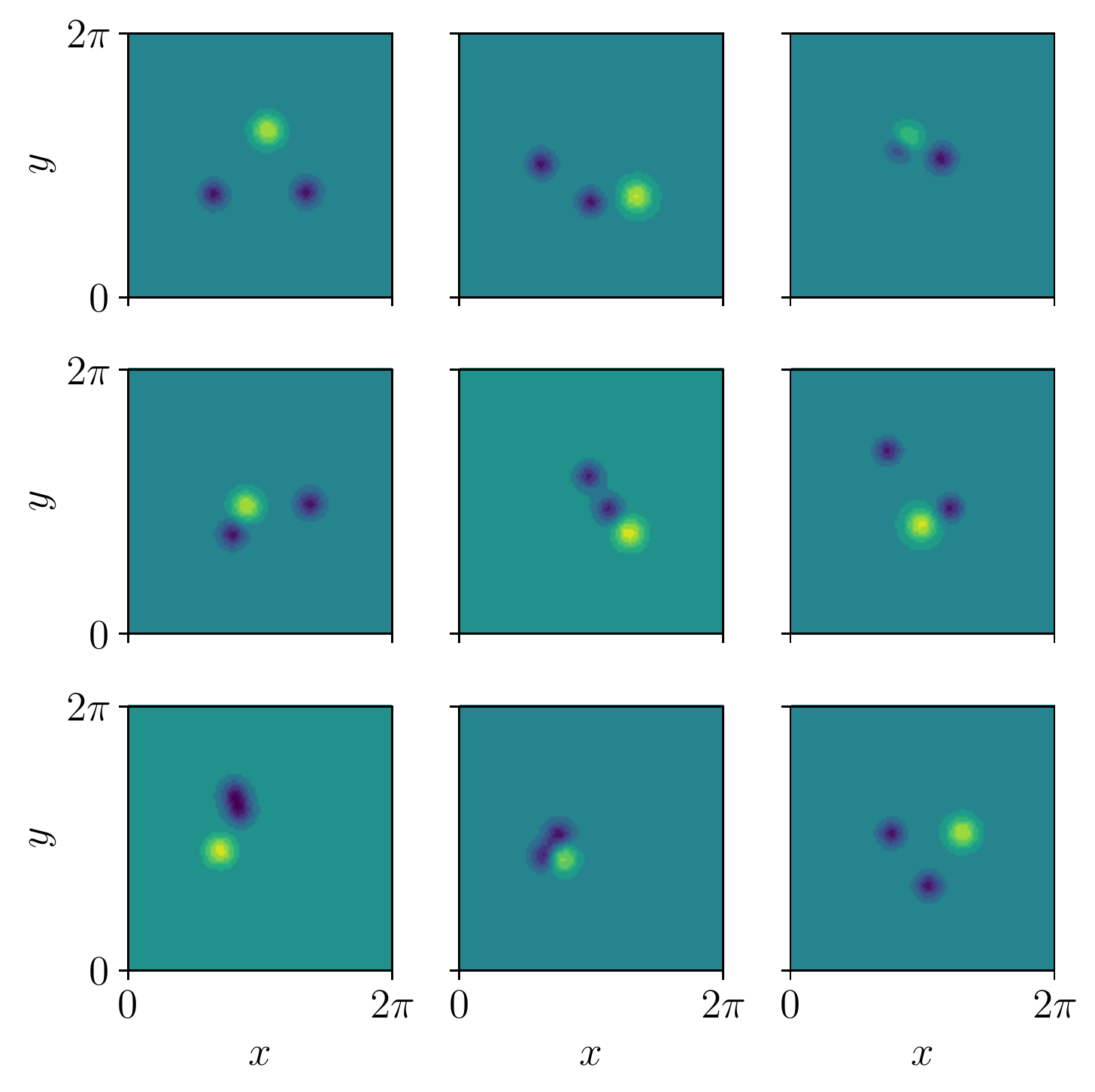}
  \caption{A set of initial conditions with three randomly located Gaussian vortices, one
  positive and two negative all with equal strength.}
  \label{fig:10}
\end{figure}

Next we will consider the problem of a two-dimensional incompressible flow in a
square periodic domain prescribed by the Navier-Stokes equations in vorticity
formulation
\begin{equation}\label{eq:30}
  \frac{\partial\omega}{\partial t} + \mathbf{u}\cdot\nabla\omega = \frac{1}{\text{Re}}\nabla^2\omega,
\end{equation}
defined on the domain $(x,y)=[0,2\pi]\times[0,2\pi]$ where $\omega(x,y,t)$ is the vorticity field,
$\mathbf{u}(x,y,t)$ is the velocity vector field. The Reynolds number is set to $\text{Re}=5\times 10^{3}$.
For the construction of the dataset, we will consider a family of initial conditions
given by a mixture of $N_v$ Gaussian vortices
\begin{equation}
  \omega(x,y,0) = \sum^{N_v}_{i=1}\delta(i)\exp\bigg(-\frac{(x-x_i)^2 + (y-y_i)^2}{0.1}\bigg),
\end{equation}
parameterized by location of the center of each vortex. Each vortex center is
sampled randomly from a square subdomain
$(x_i,y_i)\in[\pi/2,3\pi/2]\times[\pi/2,3\pi/2]\;\forall i$ as depicted in ~\autoref{fig:8}, and
the sign of each vortex is governed by $\delta(i)\in\{-1,+1\}\;\forall i$. We
consider two cases: (a) $N_v=2$, with each vortex of opposite sign, and (b)
$N_v=3$, with one positive vortex and the rest negative. A representative set of
initial conditions for the $N_v=2$ and $N_v=3$ cases can be seen in ~\autoref{fig:9}
and ~\autoref{fig:10}, respectively. This example was constructed both to showcase the application to larger scale
problems that would otherwise be too computationally intensive using a fully-connected
autoencoder and to highlight the location-invariance capabilities of the convolutional
autoencoder. The main idea is that similar to detecting an instance of an object
anywhere in an an image, the shared-weight property of each convolutional layer in
the autoencoder should be able to capture the large-parameter variations implicitly
defined in the initial condition.

To create a training dataset, ~\autoref{eq:30} is discretized pseudospectrally using
a uniform $128^2$ grid and integrated in time using the Crank-Nicholson method
to $T=250$ using a time step of $\Delta t=1\times 10^{-2}$. A parameter-varying
dataset is created by randomly sampling the initial Gaussian center locations
from a square subdomain as was previously described. Similar to first example,
after subtracting the temporal mean and feature scaling the resulting dataset
has the form
\begin{equation}\label{eq:31}
  \mathcal{X}=\{\mathbf{X}'^1_s,...,\mathbf{X}'^{N_s}_s\}\in[0,1]^{N_x\times N_y\times N_t\times N_s},
\end{equation}
where each training sample
$\mathbf{X}'^i_s=[\bm{\omega}'^1_{s,i},...,\bm{\omega}'^{N_t}_{s,i}]$ is a
matrix of two-dimensional discretized snapshots corresponding to a different set
of initial conditions. In this case, the dataset consists of a total $N_s=5120$
training samples, each with $N_t=30$ evenly sampled snapshots. Since we are
interested in employing the convolutional recurrent autoencoder, each
$\bm{\omega}'_{s,i}$ is kept as a two-dimensional array.

\begin{figure}[htbp]
  \centering
  \includegraphics[width=0.85\textwidth]{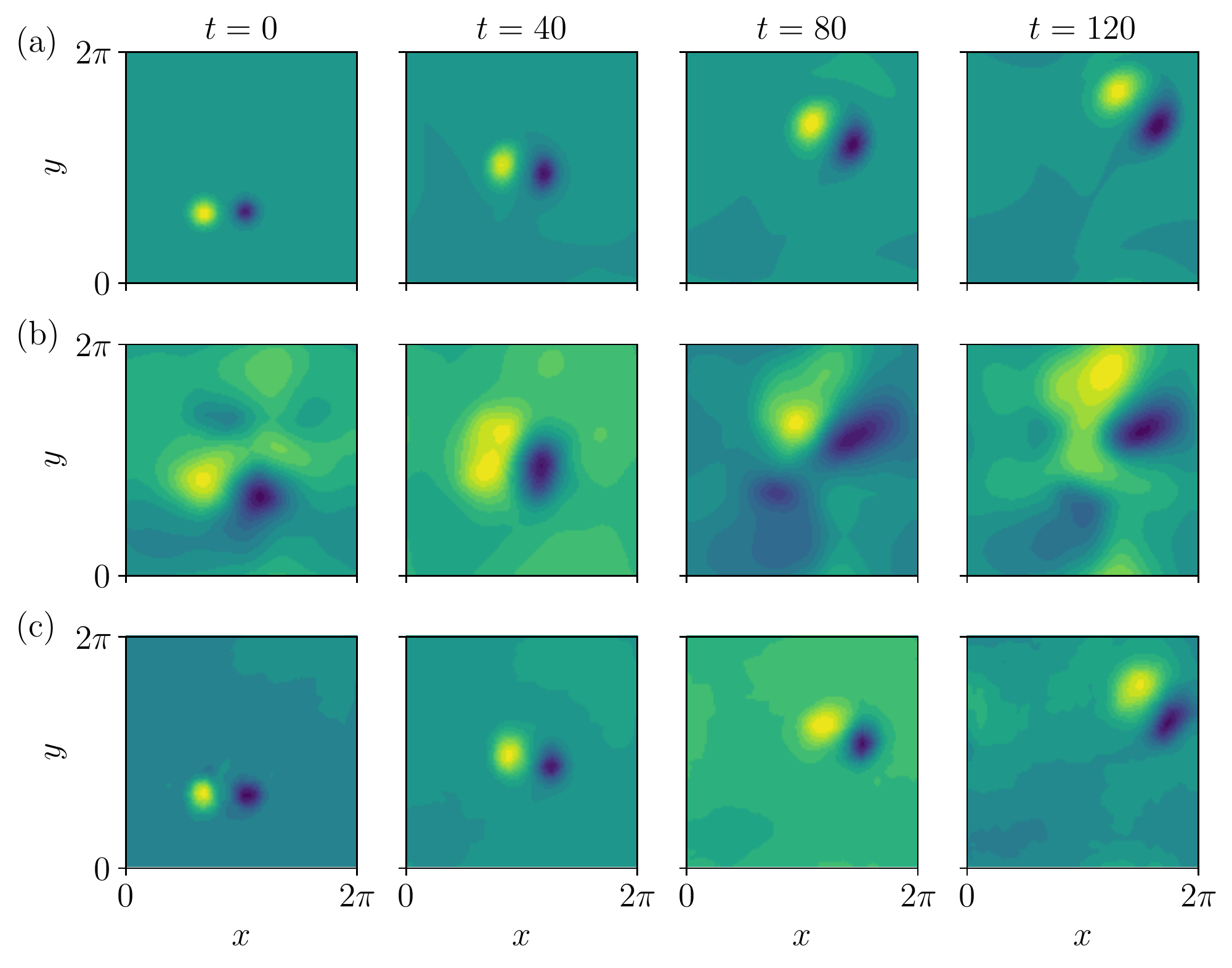}
  \caption{Comparison at $t=0,40,80,120$ of a sample trajectory using two initial vortices:
  (a) true solution, (b) rank-8 POD reconstruction using dataset with $128$
  trajectories, and (c) \textit{prediction} using a trained convolutional
  recurrent autoencoder of size $N_h=8$.}
  \label{fig:11}
\end{figure}

\begin{figure}[htbp]
  \centering
  \includegraphics[width=0.85\textwidth]{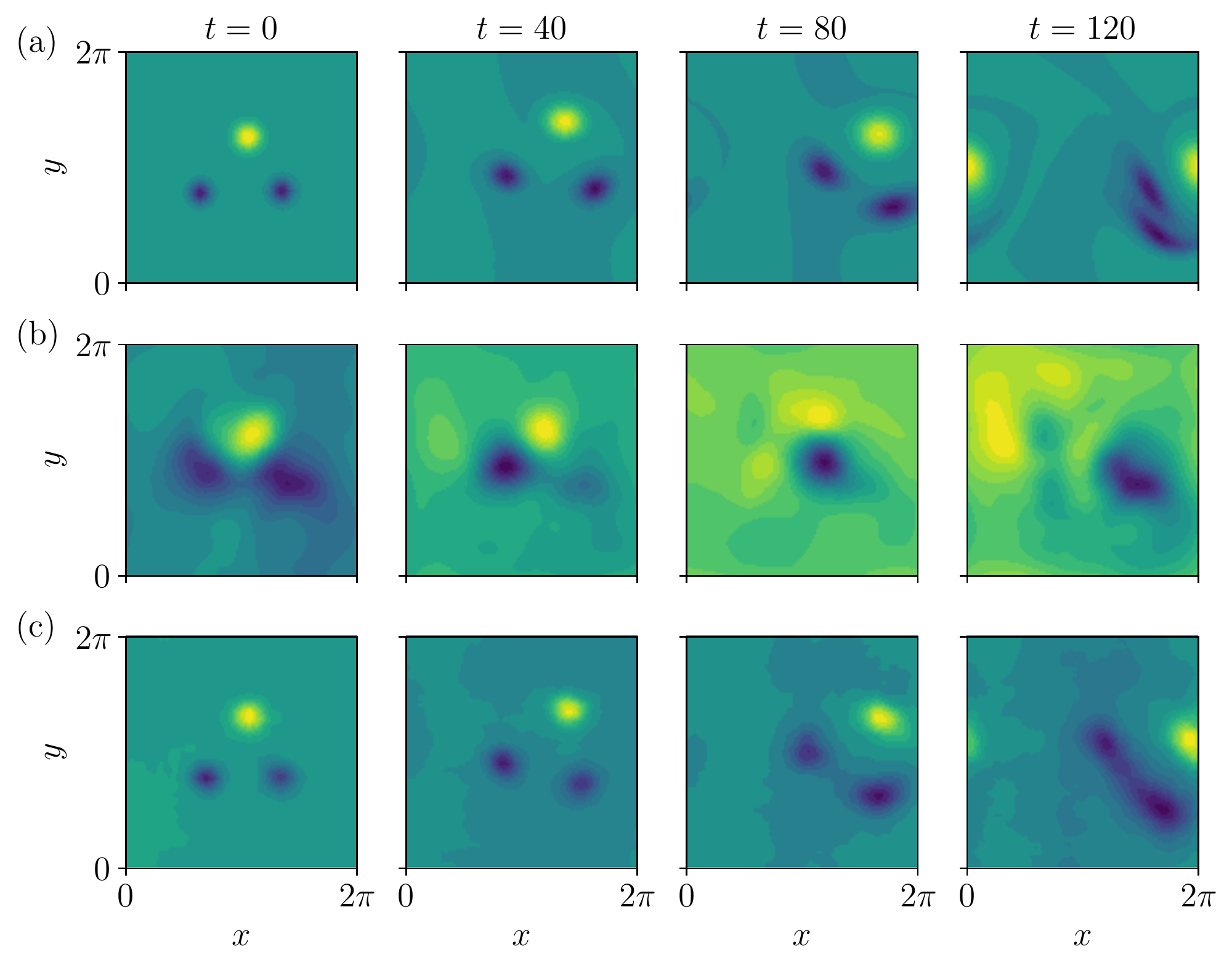}
  \caption{Comparison at $t=0,40,80,120$ of a sample trajectory using three initial vortices:
  (a) true solution, (b) rank-8 POD reconstruction using dataset with $128$
  trajectories, and (c) \textit{prediction} using a trained convolutional
  recurrent autoencoder of size $N_h=8$.}
  \label{fig:12}
\end{figure}

Three convolutional recurrent autoencoder models, with feature vector sizes
$N_h=8,16,$ and $64$, were trained using the dataset ~\autoref{eq:31} with both two
and three initial vortices. Each model was trained on an single Nvidia Tesla K20
GPU for $N_{train}=1,000,000$ iterations. Once trained, the three models were
used to predict the evolution of the vorticity field for new initial conditions.
To highlight the benefits of convolutional recurrent autoencoders for
location-invariant feature learning, we compare our prediction with a set of
best-case scenario rank-$8$ POD reconstructions.  These rank-$8$ POD
reconstructions use a dataset containing snapshots from just $128$ separate
trajectories. In this case, a rank-$8$ POD reconstruction is not sufficient to
accurately capture the correct solution since the inclusion of randomly varying
initial conditions has created a dataset that is no longer low-rank. This is
clearly shown in ~\autoref{fig:11} and ~\autoref{fig:12} for the two and three vortex
cases. The need for more and more POD modes to achieve a good reconstruction
underscores a significant disadvantage of POD-based ROMs for systems with large
variations in parameters. The convolutional recurrent autoencoder overcomes
these challenges and performs well in prediction new solutions without the need
to resort to larger-rank models. In contrast to POD-Galerkin ROMs, increasing
the number of separate trajectories in a dataset is beneficial to learning
the correct dynamic behavior.

\begin{figure}[htbp]
  \centering
  \includegraphics[width=0.5\textwidth]{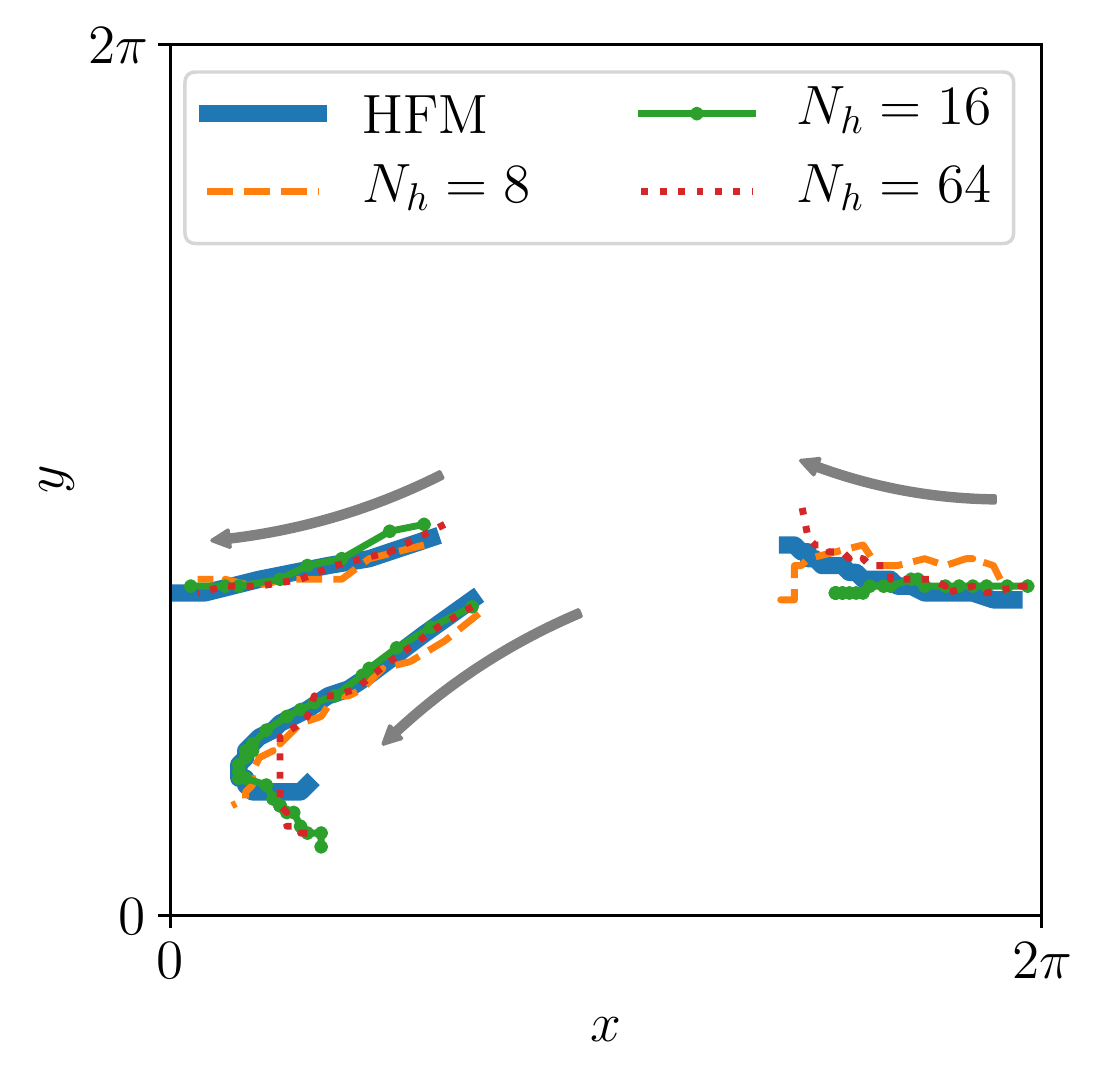}
  \caption{Evolution of the vortex centers as given the HFM solution and the predicted solutions
  using $N_h=8,16,64$.}
  \label{fig:13}
\end{figure}

\begin{figure}[htbp]
  \centering
  \includegraphics[width=0.85\textwidth]{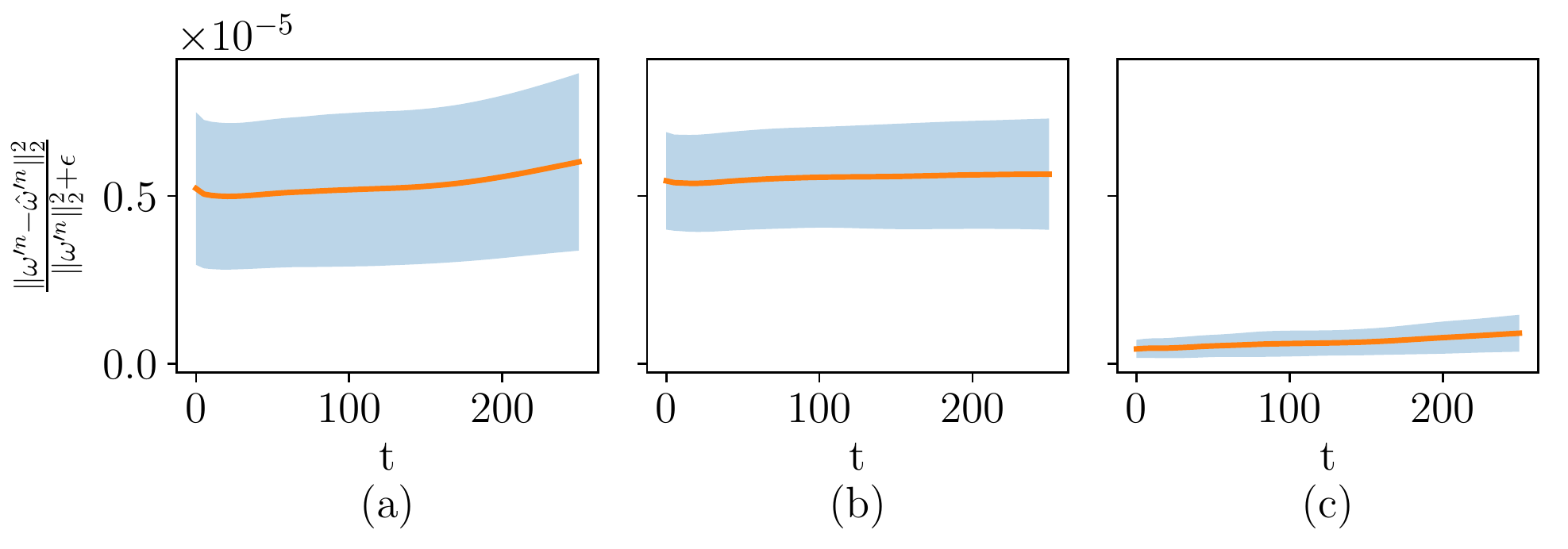}
  \caption{Mean and standard deviation of error at every time step for online
  predictions using (a) $N_h=8$, (b) $N_h=16$, and (c) $N_h=64$.}
  \label{fig:14}
\end{figure}

Similar to the first numerical example, the predictions are devoid of any
spurious oscillations that are commonplace in POD-based ROMs. Considering a
single initial condition ~\autoref{fig:13} shows the performance of each sized model
in predicting the location of each vortex as it evolves up to the training
sequence length for the two vortex case. Futher, ~\autoref{fig:14} shows the mean and
standard deviation of the scaled squared reconstruction error
\begin{equation}
  \frac{\|\bm{\omega}'^n_s-\hat{\bm{\omega}}'^n_s\|^2_F}{\|\bm{\omega}'^n_s\|^2_F+\epsilon}
\end{equation}
at each time step as calculated from 512 new prediction runs using the trained models.
In all three cases, the error did not grow significantly over the length of the
training sequence.

\subsection{Lid-driven cavity flow}
\begin{figure}[htbp]
  \centering
  \includegraphics[width=0.35\textwidth]{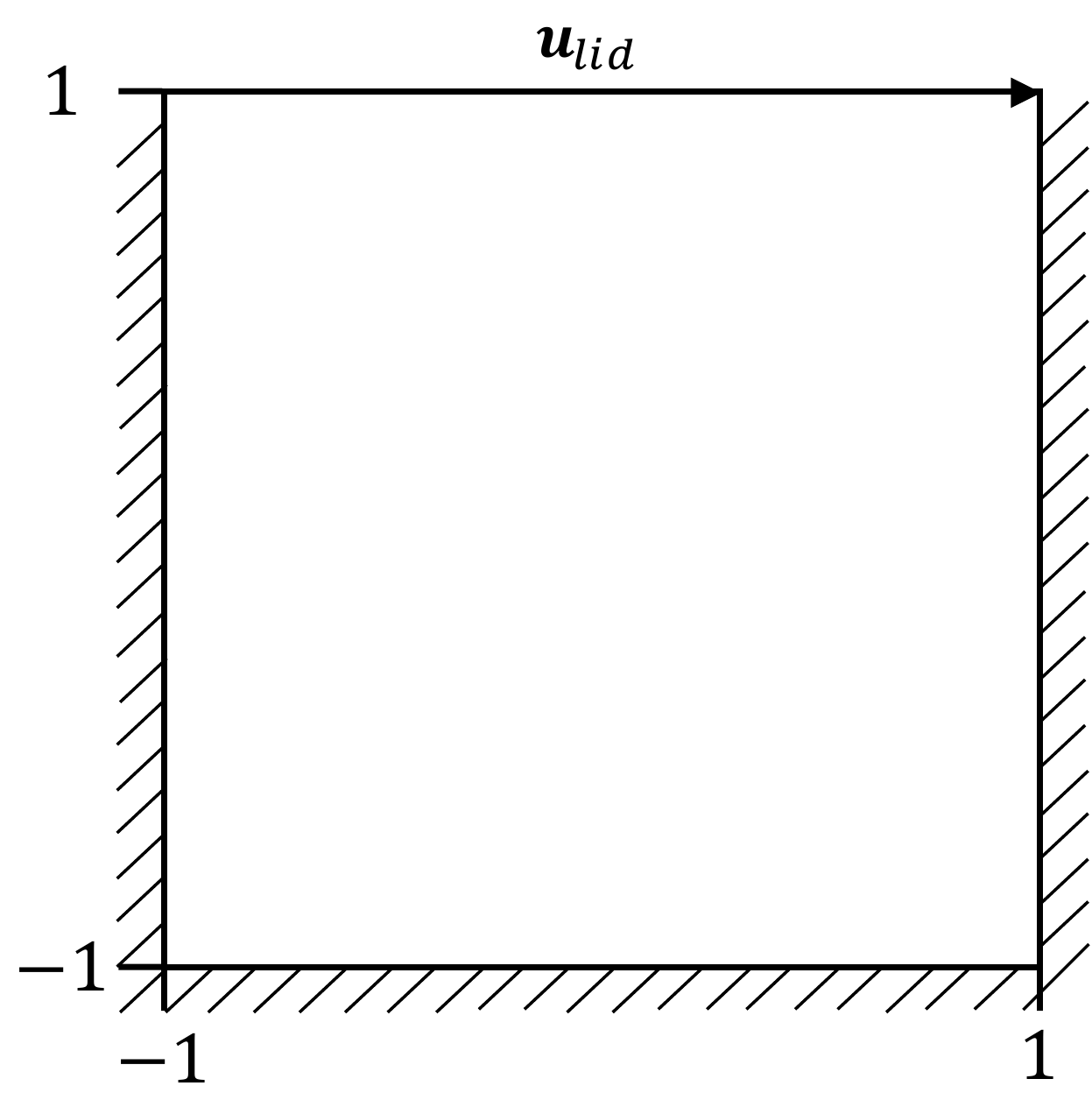}
  \caption{Lid-driven cavity domain, with lid velocity $\mathbf{u}_{\text{lid}}=(1-x^2)^2$.}
  \label{fig:15}
\end{figure}

In the final example we consider a two-dimensional incompressible flow inside a
square cavity with a lid velocity of $\mathbf{u}_{\text{lid}}=(1-x^2)^2$ at a
moderate Reynolds number $\text{Re}=2.75\times 10^4$. A graphic of the domain is
depicted in ~\autoref{fig:15}. At these Reynolds numbers the lid-driven cavity flow
is known to settle into a statistically stationary solution far from the initial
condition making it a well known benchmark for the validation of numerical
schemes and reduced order models. In particular, this benchmark is useful for
testing the stability of reduced order models ~\cite{Balajewicz2013}. The
characteristic length and velocity scales used in defining the Reynolds number
are the cavity width and the maximum lid velocity.

Consider the two-dimensional Navier-Stokes in streamfunction-vorticity formulation
defined on the square domain $(x,y)\in[-1,1]\times[-1,1]$
\begin{equation}\label{eq:32}
  \frac{\partial}{\partial t}(\nabla^2\Psi) +
  \frac{\partial\Psi}{\partial y}\frac{\partial}{\partial x}(\nabla^2\Psi) -
  \frac{\partial\Psi}{\partial y}\frac{\partial}{\partial y}(\nabla^2\Psi) = \nu\nabla^4\Psi,
\end{equation}
where $\Psi(x,y,t)$ is the streamfunction, and $\nabla^4=\nabla^2\nabla^2$ is
the biharmonic operator. To generate the training dataset, ~\autoref{eq:32} is
spatially discretized using a $128^2$ Chebyshev grid and solved numerically. The
Chebyshev coefficients are derived using the fast Fourier transform (FFT), where
the contractive nonlinearities are handled pseudospectrally. The equations are
integrated in time using a semi-implicit, second order Euler scheme. Since the
statistically stationary solution is far from the initial condition, we first
initialize the simulation over $7,500,000$ time steps with time-step size
$\Delta t = 1\times 10^{-4}$. The following $2,500,000$ time steps are then used
to create a dataset in the form of ~\autoref{eq:26} with $N_s=1110$ where now each
training sample is
\begin{equation}\label{eq:33}
  \mathbf{X}'^i_s=[\bm{\Psi}'^i_s,\bm{\Psi}'^{i+m}_s,\bm{\Psi}'^{i+2m}_s,...,\bm{\Psi}'^{i+(N_t-1)m}_s]\in\mathbb{R}^{N_x\times N_y\times N_t},
\end{equation}
where each $\bm{\Psi}'^i_s$ is a discretized two-dimensional snapshot of
~\autoref{eq:32}, $N_t=35$, and $m$ is taken to be $100$. In doing this, we ensure
that the initial training snapshot used to initialize the RNN portion of the
model evenly samples the entire trajectory of ~\autoref{eq:32}. The result is the
construction a training dataset that gives a good representation of the dynamics
for the RNN to learn. In addition, an interpolation step onto a uniform $128^2$
is performed to ensure each filter $\mathbf{K}^f$ acts on equally
physically-sized receptive fields.

\begin{figure}[htbp]
  \centering
  \includegraphics[width=0.85\textwidth]{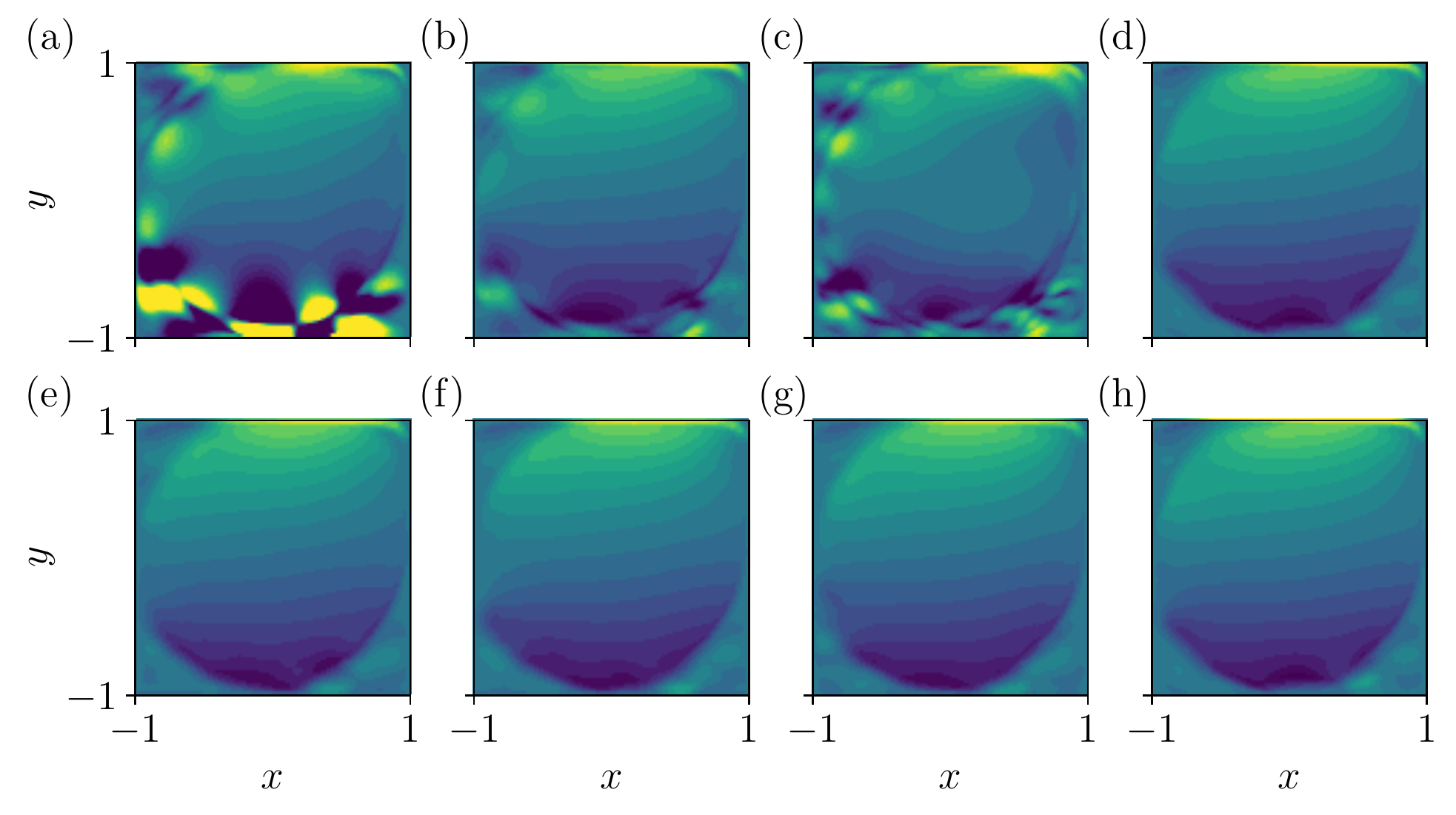}
  \caption{$u(x,y,t)$ contours of the lid-driven cavity flow at $t=250\;s$ using
  the optimal POD reconstruction: (a) $N_h=8$ (note: $t=60$ shown, right before blowup), (b) $N_h=16$, (c) $N_h=64$,
  (d) true solution; and \textit{predicted} contours using the
  convolutional recurrent autoencoder model with hidden state sizes
  (e) $N_h=8$, (f) $N_h=16$, (g) $N_h=64$, (h) true solution.}
  \label{fig:16}
\end{figure}

\begin{figure}[htbp]
  \centering
  \includegraphics[width=0.85\textwidth]{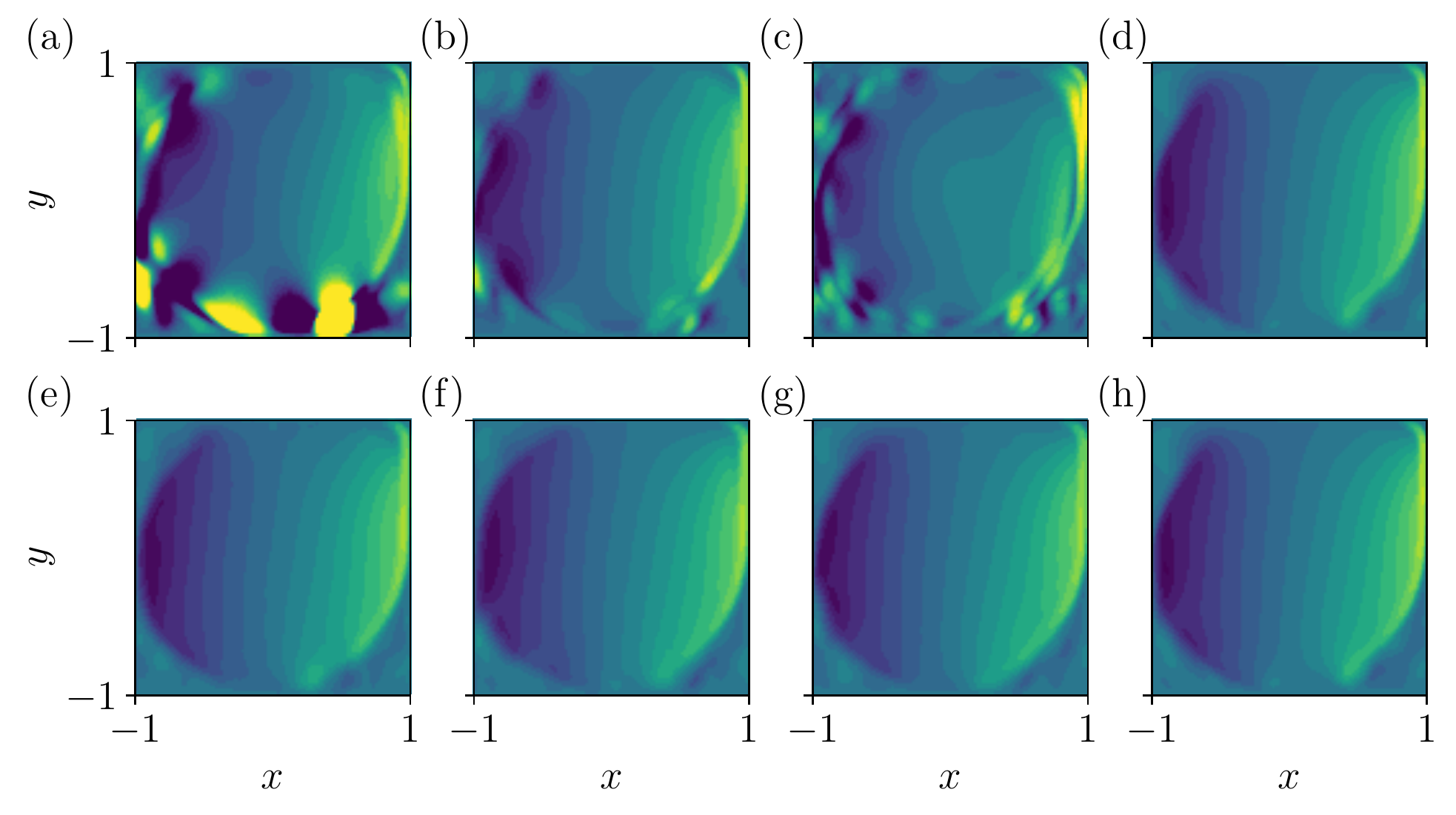}
  \caption{$v(x,y,t)$ contours of the lid-driven cavity flow at $t=250\;s$ using
  the optimal POD reconstruction: (a) $N_h=8$ (note: $t=60$ shown, right before blowup), (b) $N_h=16$, (c) $N_h=64$,
  (d) true solution; and \textit{predicted} contours using the
  convolutional recurrent autoencoder model with hidden state sizes
  (e) $N_h=8$, (f) $N_h=16$, (g) $N_h=64$, (h) true solution.}
  \label{fig:17}
\end{figure}

\begin{figure}[htbp]
  \centering
  \includegraphics[width=0.85\textwidth]{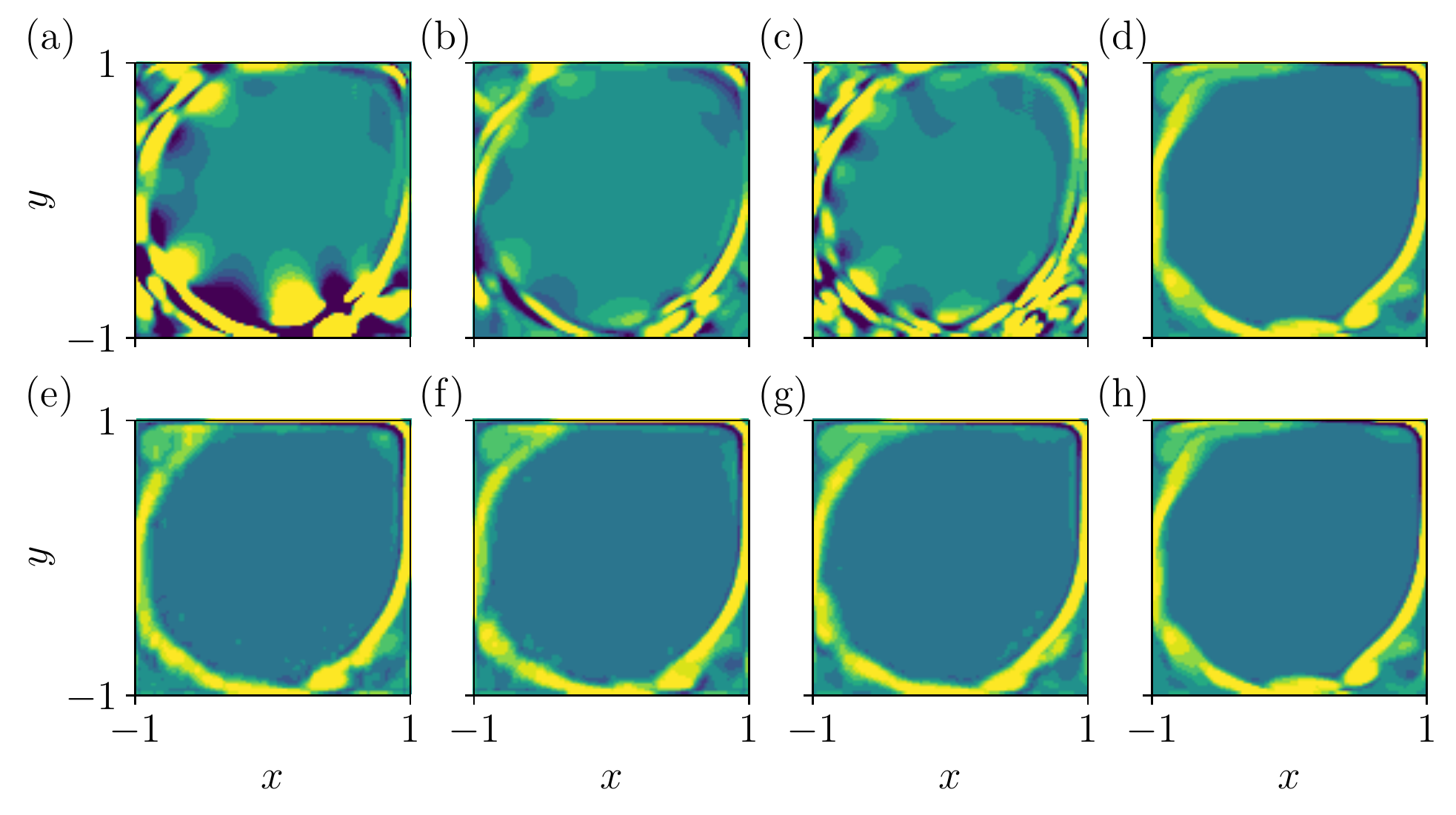}
  \caption{Vorticity contours of the lid-driven cavity flow at $t=250\;s$ using
  the optimal POD reconstruction: (a) $N_h=8$ (note: $t=60$ shown, right before blowup), (b) $N_h=16$, (c) $N_h=64$,
  (d) true solution; and \textit{predicted} contours using the
  convolutional recurrent autoencoder model with hidden state sizes
  (e) $N_h=8$, (f) $N_h=16$, (g) $N_h=64$, (h) true solution.}
  \label{fig:18}
\end{figure}

Three convolutional recurrent autoencoder models were trained using this dataset,
again with low-dimensional representations of sizes $N_h=8,16,$ and $64$. In this
case all three models were trained on a single Nvidia Tesla K20 GPU for $N_{train}=600,000$
iterations. The online performance of the each model was evaluated by initializing
each model with a slightly perturbed version of the first snapshot of the entire
dataset and evaluating for 2500 prediction steps, over 70 times the length
of each training sequence. We perform the same with three equivalently sized POD-Galerkin ROMs.
~\autoref{fig:16,fig:17,fig:18} depict the final predicted velocity fields $u(x,y,t)$, $v(x,y,t)$, as well
as the predicted vorticity field $\omega(x,y,t)$ using traditional POD-Galerkin ROMs
and our convolutional recurrent autoencoder model for $N_h=8,16,$ and $64$. In
all three reconstructed fields the poor performance of POD-Galerkin ROMs can be
easily noticed by the spurious oscillations present in the field. This is
in contrast to the predictions presented using our approach, which nearly capture
the exact solution even after long-term prediction.

\begin{figure}[htbp]
  \centering
  \includegraphics[width=0.85\textwidth]{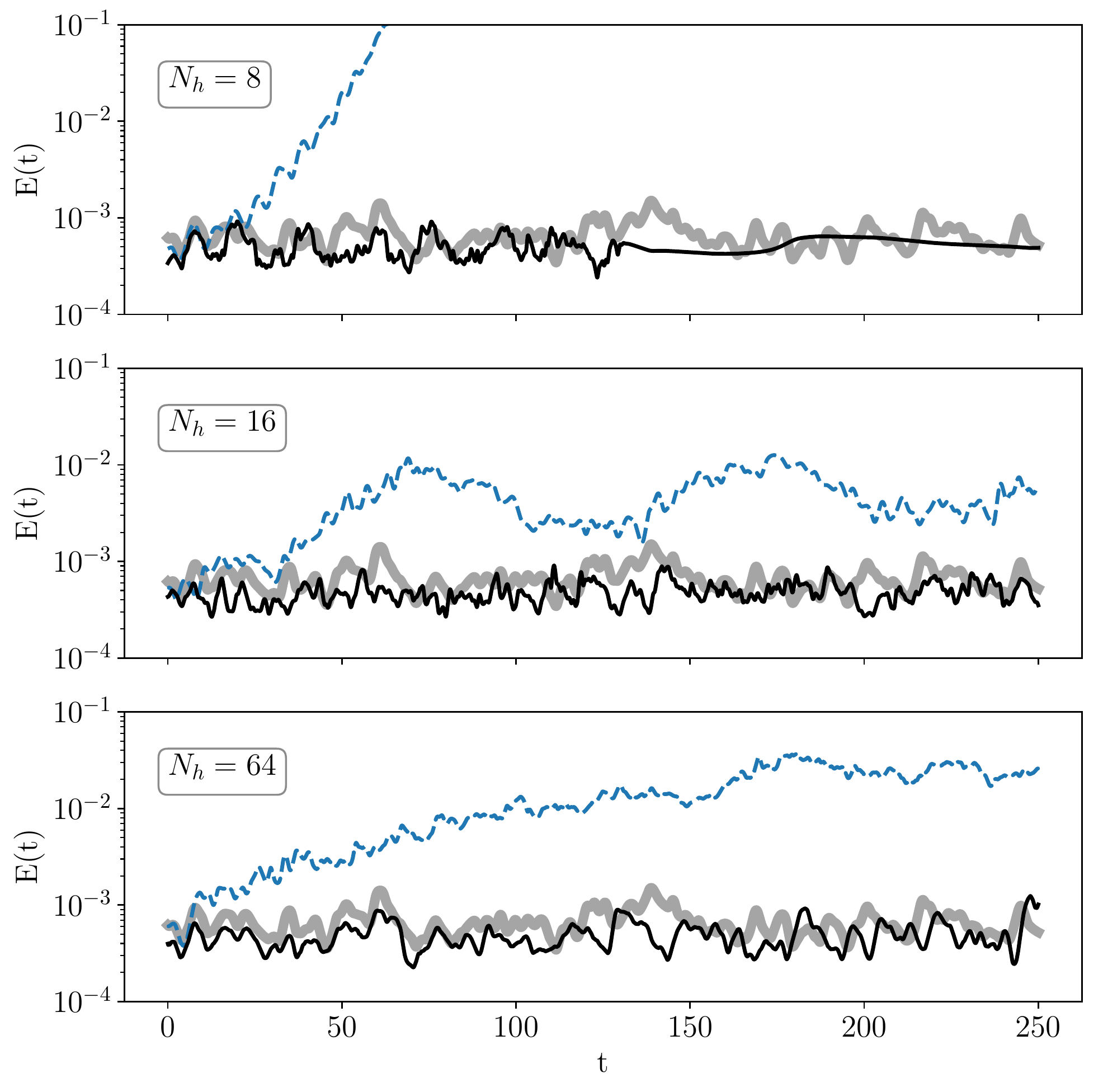}
  \caption{The evolution of the instantaneous turbulent kinetic energy for the
  lid-driven cavity flow from the DNS (thick grey lines), standard POD-based
  Galerkin ROMs (blue dashed lines), and our method (solid black lines).}
  \label{fig:19}
\end{figure}

In fact, we only present predictions up until $t=60$ for the $N_h=8$ POD-Galerkin ROM since
instabilities cause the solution to diverge. This can be seen more clearly
in ~\autoref{fig:19}, which compares the instantaneous turbulent kinetic energy (TKE) of the flow
\begin{equation}\label{eq:34}
  E(t) = \frac{1}{2}\int_{\Omega}\bigg(u(t)'^2 + v(t)'^2\bigg)d\Omega
\end{equation}
where $u(t)'$ and $v(t)'$ are the instantaneous velocity fluctuations around the
mean and $\Omega$ represents the fluid domain. The TKE can be seen as a measure
of the energy content within the flow. For statistically stationary flows, such
as the one considered in this example, the TKE should hover around a mean value.
In ~\autoref{fig:19} we see that the POD-Galerkin models fail to capture the correct
TKE, and in the case of $N_h=8$ instabilities lead to eventual divergence.

\begin{figure}[htbp]
  \centering
  \includegraphics[width=0.7\textwidth]{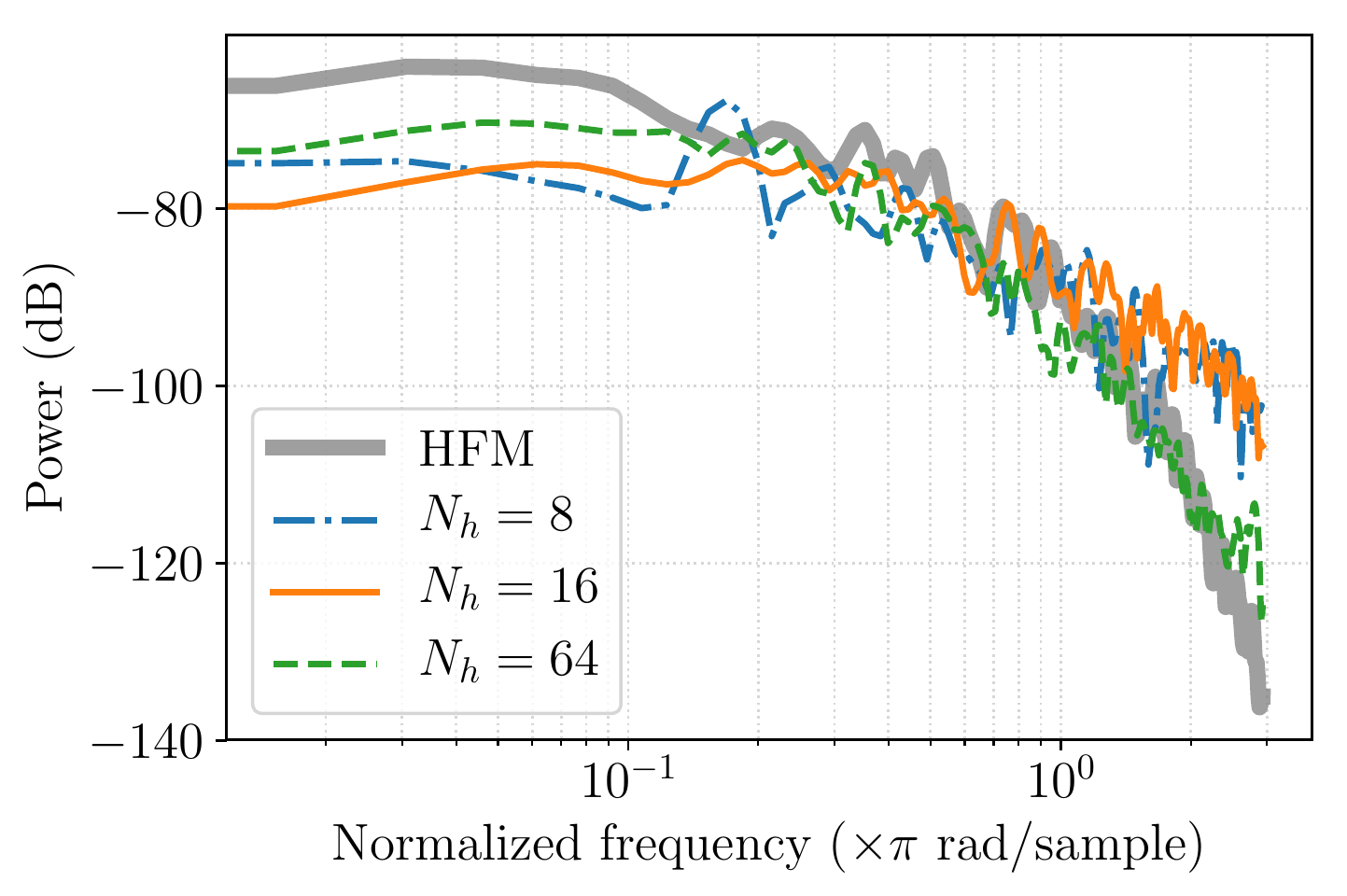}
  \caption{PSD of the turbulent kinetic energy of the lid-driven cavity flow.}
  \label{fig:20}
\end{figure}

Against this backdrop, we can see that our approach vastly outperforms traditional
POD-Galerkin ROMs. All velocity and vorticity reconstructions are in good agreement
with the HFM solution. As the size of the model increases to $N_h=64$, we see that predicted TKE
is in good agreement with that of the HFM. It should be noted that the lid-driven
cavity flow at these Reynolds numbers exhibits chaotic motion, thus a best-case
scenario would be to capture the right TKE in a statistical sense. This can be
seen further in ~\autoref{fig:20} which compares the power spectral density of each
predicted TKE with that of the HFM.

While each model prediction capatures the general behavior of the HFM, there is
some high spatial frequency error evident throughout the domain in each reconstruction.
Interestingly, the stability of the RNN portion of the each model remains unaffected
by this high-frequency noise suggesting that it is due only to the transpose
convolutional decoder. This is possibly a result of performing a strided transpose
convolution at each layer of the decoder. It is possible and perhaps beneficial
to include a final undilated convolutional layer with a single feature map to filter
some of the high-frequency reconstruction noise.

\section{Conclusions}
\label{sec:conclusions}
In this work we propose a completely data-driven nonlinear reduced order model
based on a convolutional recurrent autoencoder architecture for application to
parameter-varying systems and systems requiring long-term stability. The
construction of the convolutional recurrent autoencoder consists of two major
components each of which performs a key task in projection based reduced
order modeling. First a convolutional autoencoder is designed to identify a
low-dimensional representation of two-dimensional input data in terms of
intrinsic coordinates on some low-dimensional manifold embedded in the original,
high-dimensional space. This is done by considering a 4-layer convolutional
encoder which computes a hierarchy of localized, location invariant features
that are passed to a two-layer fully connected encoder. The result of this is a
mapping from the high-dimensional input space to a low-dimensional
data-supporting manifold. An equivalent decoder architecture is considered for
efficiently mapping from the low-dimensional representation to the original
space. This can be intuitively understood as a nonlinear generalization of POD,
where the structure of the manifold is more expressive than the linear subspaces
learned by POD-based methods. The second important component of the proposed
convolutional recurrent autoencoder is a modified version of an LSTM network
which models the dynamics on the manifold learned by the autoencoder. The LSTM
network is modified to require only information from the low-dimensional
representation thereby avoiding costly reconstruction of the full state at every
evolution step.

An offline training and online prediction strategy for the convolutional
recurrent autoencoder is proposed in this work. The training algorithm exploits
the modularity of the model by splitting each forward pass into two steps. The
first step running a forward pass on the autoencoder while creating a temporary
batch of target low-dimensional representations which are then used in the second
step, which is the forward pass of the modified LSTM network. The backwards pass,
or parameter update is then performed jointly equally weighting autoencoder
reconstruction error and the prediction error of the modified LSTM network.

We demonstrated our approach on three illustrative nonlinear model reduction
examples. The first emphasizes the expressive power of using fully-connected
autoencoders equipped with nonlinear activation functions on performing model
reduction tasks in contrast to POD-based methods. The second highlights the
performance of the convolutional recurrent autoencoder, and in particular its
location-invariant properties, in parametric model reduction with initial
condition exhibiting large parameter variations. The final example demonstrates
the stability of convolutional recurrent autoencoders when performing long-term
predictions of choatic incompressible flows. Collectively, these numerical
examples show that our convolutional recurrent autoencoder model outperforms
traditional POD-Galerkin ROMs both in terms of prediction quality, parameter
variations, and stability while also offering other advantages such as location
invariant feature learning and non-intrusiveness. In fact, although in this work
we make use of canonical model reduction examples based on computational physics
problems, our approach is completely general and can be applied to arbitrary
high-dimensional spatiotemporal data. When compared to existing
autoencoder-based reduced order modeling strategies, our model provides access
to larger-sized problems while keeping the number of trainable parameters low
compared to fully-connected autoencoders.

\subsection{Future work}

This work shows the feasibility of using deep learning-based strategies for
performing nonlinear model reduction and more generally modeling complex
dynamical system in a completely data-driven and non-intrusive manner. Although
this work presents promising predictive results for both parameter-varying model
reduction problems and problems requiring long-term stability, these methods
remain in their infancy and their full capabilities are yet unknown. There are
multiple directions in which this work can be extended. One such direction is
improving the design of the convolutional transpose decoder. As it stands, the
main source of error in our results is high-frequency in nature an appears only
during the decoding phase. Considering this, future decoder designs could
include more efficient filtering strategies. Another possible direction is in
the dynamic modeling of the low-dimensional representations. In this work, we
considered samples with spatial parameter variations and thus the design of the
LSTM network could remain unchanged. However, there is potential for deep
learning-based dynamic modeling approaches that exploit multi-scale phenomena
inherent in many physical systems. Finally, a much more challenging problem is
the reconciliation of deep learning-based performance gains with physical
intuition. This issue permeates throughout all fields where deep learning has
made an impact: what is it actually doing? Developing our understanding of deep
learning-based modeling strategies can potentially provide us with deeper
insights of the dynamics inherent in a physical system.

\section*{Acknowledgments}
This material is based upon the work supported by the Air Force Office of
Scientific Research under Grant No. FA9550-17-1-0203. Simulations and model
training were also made possible in part by an exploratory award from the Blue
Waters sustained-petascale computing project, which is supported by the National
Science Foundation (awards OCI-0725070 and ACI-1238993) and the state of
Illinois. Blue Waters is a joint effort of the University of Illinois at
Urbana-Champaign and its National Center for Supercomputing Applications.

\end{document}